
\overfullrule=0pt

\def\bkt{\it} 

\font\headerfont=cmr10
\font\tenbbold=msbm10

\font\sectionfont=cmbx12

\newlinechar=`@
\def\forwardmsg#1#2#3{\immediate\write16{@*!*!*!* forward reference should
be: @\noexpand\forward{#1}{#2}{#3}@}}
\def\nodefmsg#1{\immediate\write16{@*!*!*!* #1 is an undefined reference@}}

\def\forwardsub#1#2{\def\newref{{#2}{#1}}}

\def\forward#1#2#3{%
\expandafter\expandafter\expandafter\forwardsub\expandafter{#3}{#2}
\expandafter\ifx\csname#1\endcsname\relax\else%
\expandafter\ifx\csname#1\endcsname\newref\else%
\expandafter\ifx\csname#2\endcsname\relax\else%
\forwardmsg{#1}{#2}{#3}\fi\fi\fi%
\expandafter\let\csname#1\endcsname\newref}

\def\firstarg#1{\expandafter\argone #1}\def\argone#1#2{#1}
\def\secondarg#1{\expandafter\argtwo #1}\def\argtwo#1#2{#2}

\def\ref#1{\expandafter\ifx\csname#1\endcsname\relax%
  {\nodefmsg{#1}\bf`#1'}\else%
  \expandafter\firstarg\csname#1\endcsname%
  ~\htmllocref{#1}{\expandafter\secondarg\csname#1\endcsname}\fi}

\def\refscor#1{\expandafter\ifx\csname#1\endcsname\relax
  {\nodefmsg{#1}\bf`#1'}\else
  Corollaries~\htmllocref{#1}{\expandafter\secondarg\csname#1\endcsname}\fi}

\def\refs#1{\expandafter\ifx\csname#1\endcsname\relax
  {\nodefmsg{#1}\bf`#1'}\else
  \expandafter\firstarg\csname #1\endcsname
  s~\htmllocref{#1}{\expandafter\secondarg\csname#1\endcsname}\fi}

\def\refn#1{\expandafter\ifx\csname#1\endcsname\relax
  {\nodefmsg{#1}\bf`#1'}\else
  \htmllocref{#1}{\expandafter\secondarg\csname #1\endcsname}\fi}

\def\pageref#1{\expandafter\ifx\csname#1\endcsname\relax
  {\nodefmsg{#1}\bf`#1'}\else
  \expandafter\firstarg\csname#1\endcsname%
  ~\htmllocref{#1}{\expandafter\secondarg\csname#1\endcsname}\fi}

\def\pagerefs#1{\expandafter\ifx\csname#1\endcsname\relax
  {\nodefmsg{#1}\bf`#1'}\else
  \expandafter\firstarg\csname#1\endcsname%
  s~\htmllocref{#1}{\expandafter\secondarg\csname#1\endcsname}\fi}

\def\pagerefn#1{\expandafter\ifx\csname#1\endcsname\relax
  {\nodefmsg{#1}\bf`#1'}\else
  \htmllocref{#1}{\expandafter\secondarg\csname#1\endcsname}\fi}






\newif\ifhypers 
\hypersfalse


\edef\freehash{\catcode`\noexpand\#=\the\catcode`\#}%
\catcode`\#=12
\freehash
\let\freehash=\relax
\ifhypers\fi
\def\puthtml#1{\ifhypers\fi}
\def\htmlanchor#1#2{\puthtml{<a name="#1">}#2\puthtml{</a>}}
\def\@pdfm@mark#1{}
\def\setlink#1{\colored{\linkcolor}{#1}}%
\def\setlink#1{\ifdraft\Purple{#1}\else{#1}\fi}%

%

%
\def\htmllocref#1#2{\ifhypers\leavevmode\fi\setlink{#2}\ifhypers\fi\relax}%
%
%
%
%
\def\Acrobatmenu#1#2{%
  \@pdfm@mark{%
    bann <<
      /Type /Annot
      /Subtype /Link
      /A <<
        /S /Named
        /N /#1
      >>
      /Border [\@pdfborder]
      /C [\@menubordercolor]
    >>%
   }%
  \Hy@colorlink{\@menucolor}#2\Hy@endcolorlink
  \@pdfm@mark{eann}%
}
\def\@pdfborder{0 0 1}
\def\@menubordercolor{1 0 0}
\def\@menucolor{red}

\def\ifempty#1#2\endB{\ifx#1\endA}
\def\makeref#1#2#3{\ifempty#1\endA\endB\else\forward{#1}{#2}{#3}\fi\unskip}


\def\crosshairs#1#2{
  \dimen1=.002\drawx
  \dimen2=.002\drawy
  \ifdim\dimen1<\dimen2\dimen3\dimen1\else\dimen3\dimen2\fi
  \setbox1=\vclap{\vline{2\dimen3}}
  \setbox2=\clap{\hline{2\dimen3}}
  \setbox3=\hstutter{\kern\dimen1\box1}{4}
  \setbox4=\vstutter{\kern\dimen2\box2}{4}
  \setbox1=\vclap{\vline{4\dimen3}}
  \setbox2=\clap{\hline{4\dimen3}}
  \setbox5=\clap{\copy1\hstutter{\box3\kern\dimen1\box1}{6}}
  \setbox6=\vclap{\copy2\vstutter{\box4\kern\dimen2\box2}{6}}
  \setbox1=\vbox{\offinterlineskip\box5\box6}
  \smash{\vbox to #2{\hbox to #1{\hss\box1}\vss}}}

\def\boxgrid#1{\rlap{\vbox{\offinterlineskip
  \setbox0=\hline{\wd#1}
  \setbox1=\vline{\ht#1}
  \smash{\vbox to \ht#1{\offinterlineskip\copy0\vfil\box0}}
  \smash{\vbox{\hbox to \wd#1{\copy1\hfil\box1}}}}}}

\def\drawgrid#1#2{\vbox{\offinterlineskip
  \dimen0=\drawx
  \dimen1=\drawy
  \divide\dimen0 by 10
  \divide\dimen1 by 10
  \setbox0=\hline\drawx
  \setbox1=\vline\drawy
  \smash{\vbox{\offinterlineskip
    \copy0\vstutter{\kern\dimen1\box0}{10}}}
  \smash{\hbox{\copy1\hstutter{\kern\dimen0\box1}{10}}}}}

\long\def\boxit#1#2#3{\hbox{\vrule
  \vtop{%
    \vbox{\hsize=#2\hrule\kern#1%
      \hbox{\kern#1#3\kern#1}}%
    \kern#1\hrule}%
    \vrule}}
\long\def\boxitv#1#2#3{\boxit{#1}{#2}{{\hbox to #2{\vrule%
  \vbox{\hsize=#2\hrule\kern#1%
      \hbox{\kern#1#3\kern#1}}%
    \kern#1\hrule}%
    \vrule}}}

\newdimen\boxingdimen
\long\def\boxing#1{\boxingdimen=\hsize\advance\boxingdimen by -2ex%
\vskip0.5cm\boxit{1ex}{\boxingdimen}{\vbox{#1}}\vskip0.5cm}

\newdimen\boxrulethickness \boxrulethickness=.4pt
\newdimen\boxedhsize
\newbox\textbox
\newdimen\originalbaseline
\newdimen\hborderspace\newdimen\vborderspace
\hborderspace=3pt \vborderspace=3pt

\def\preparerulebox#1#2{\setbox\textbox=\hbox{#2}%
   \originalbaseline=\vborderspace
   \advance\originalbaseline\boxrulethickness
   \advance\originalbaseline\dp\textbox
   \def\Borderbox{\vbox{\hrule height\boxrulethickness
     \hbox{\vrule width\boxrulethickness\hskip\hborderspace
     \vbox{\vskip\vborderspace\relax#1{#2}\vskip\vborderspace}%
     \hskip\hborderspace\vrule width\boxrulethickness}%
     \hrule height\boxrulethickness}}}

\def\hrulebox#1{\hbox{\preparerulebox{\hbox}{#1}%
   \lower\originalbaseline\Borderbox}}
\def\vrulebox#1#2{\vbox{\preparerulebox{\vbox}{\hsize#1#2}\Borderbox}}
\def\parrulebox#1\par{\boxedhsize=\hsize\advance\boxedhsize
   -2\boxrulethickness \advance\boxedhsize-2\hborderspace
   \vskip3\parskip\vrulebox{\boxedhsize}{#1}\vskip\parskip\par}
\def\sparrulebox#1\par{\vskip1truecm\boxedhsize=\hsize\advance\boxedhsize
   -2\boxrulethickness \advance\boxedhsize-2\hborderspace
   \vskip3\parskip\vrulebox{\boxedhsize}{#1}\vskip.5truecm\par}



%

\newcount\sectno \sectno=0
\newcount\appno \appno=64 
\newcount\thmno \thmno=0
\newcount\randomct \newcount\rrandomct
\newif\ifsect \secttrue
\newif\ifapp \appfalse
\newif\ifnumct \numcttrue
\newif\ifchap \chapfalse
\newif\ifdraft \draftfalse
\newif\ifcolor\colorfalse
\ifcolor\input colordvi\else\input blackdvi\fi%


\def\widow#1{\vskip 0pt plus#1\vsize\goodbreak\vskip 0pt plus-#1\vsize}


\def\stdskip{\vskip 9pt plus3pt minus 3pt}
\def\stdbreak{\par\removelastskip\penalty-100\stdskip}
\def\halfbreak{\vskip 0.6ex\penalty-100}

\def\proof{\stdbreak\noindent%
  \ifdraft\let\proofcolor=\Purple\else\let\proofcolor=\Black\fi%
  \proofcolor{{\sl Proof. }}}

\def\proofone{\stdbreak\noindent%
  \ifdraft\let\proofcolor=\Purple\else\let\proofcolor=\Black\fi%
  \proofcolor{{\sl Proof \#1:\ \ }}}
\def\prooftwo{\stdbreak\noindent%
  \ifdraft\let\proofcolor=\Purple\else\let\proofcolor=\Black\fi%
  \proofcolor{{\sl Proof \#2:\ \ }}}
\def\proofthree{\stdbreak\noindent%
  \ifdraft\let\proofcolor=\Purple\else\let\proofcolor=\Black\fi%
  \proofcolor{{\sl Proof \#3:\ \ }}}
\def\proofof#1{\stdbreak\noindent%
  \ifdraft\let\proofcolor=\Purple\else\let\proofcolor=\Black\fi%
  \proofcolor{{\sl Proof of #1:\ }}}
\def\claim{\stdbreak\noindent%
  \ifdraft\let\proofcolor=\Purple\else\let\proofcolor=\Black\fi%
  \proofcolor{{\sl Claim:\ }}}
\def\proofclaim{\stdbreak\noindent%
  \ifdraft\let\proofcolor=\Purple\else\let\proofcolor=\Black\fi%
  \proofcolor{{\sl Proof of the claim:\ }}}

\newif\ifnumberbibl \numberbibltrue
\newdimen\labelwidth

\def\references{\bgroup
  \edef\numtoks{}%
  \global\thmno=0
  \setbox1=\hbox{[999]} 
  \labelwidth=\wd1 \advance\labelwidth by 2.5em
  \ifnumberbibl\advance\labelwidth by -2em\fi
  \parindent=\labelwidth \advance\parindent by 0.5em 
  \removelastskip
  \widow{0.1}
  \vskip 24pt plus 6pt minus 6 pt
  \frenchspacing
  \immediate\write\isauxout{\noexpand\forward{Bibliography}{}{\the\pageno}}%
  \immediate\write\iscontout{\noexpand\contnosectlist{}{Bibliography}{\the\pageno}}%
  \ifdraft\let\refcolor=\Maroon\else\let\refcolor=\Black\fi%
  \leftline{\sectionfont\refcolor{References}}
  \ifhypers%
     \global\thmno=0\relax%
     \ifsect%
	\global\advance\sectno by 1%
	\edef\numtoks{\number\sectno}%
  	\hbox{}%
     \else%
	\edef\numtoks{}%
  	\hbox{}%
     \fi%
  \fi%
  \nobreak\stdskip\noindent}%
\def\endreferences{\nonfrenchspacing\egroup}

\def\referencesn{\bgroup
  \edef\numtoks{}%
  \global\thmno=0
  \setbox1=\hbox{[999]} 
  \labelwidth=\wd1 \advance\labelwidth by 2.5em
  \ifnumberbibl\advance\labelwidth by -2em\fi
  \parindent=\labelwidth \advance\parindent by 0.5em 
  \removelastskip
  \widow{0.1}
  \vskip 24pt plus 6pt minus 6 pt
  \frenchspacing
  \immediate\write\isauxout{\noexpand\forward{Bibliography}{}{\the\pageno}}%
  \immediate\write\iscontout{\noexpand\contnosectlist{}{Bibliography}{\the\pageno}}%
  \ifdraft\let\refcolor=\Maroon\else\let\refcolor=\Black\fi%
  \ifhypers%
     \global\thmno=0\relax%
     \ifsect%
	\global\advance\sectno by 1%
	\edef\numtoks{\number\sectno}%
  	\hbox{}%
     \else%
	\edef\numtoks{}%
  	\hbox{}%
     \fi%
  \fi%
  \nobreak\stdskip\noindent}%

\def\bitem#1{\global\advance\thmno by 1%
  \ifdraft\let\itemcolor=\Purple\else\let\itemcolor=\Black\fi%
  \outer\expandafter\def\csname#1\endcsname{\the\thmno}%
  \edef\numtoks{\number\thmno}%
  \ifhypers\htmlanchor{#1}{\makeref{\noexpand#1}{REF}{\numtoks}}\fi%
  \ifnumberbibl
    \immediate\write\isauxout{\noexpand\forward{\noexpand#1}{}{\the\thmno}}
    \item{\hbox to \labelwidth{\itemcolor{\hfil\the\thmno.\ \ }}}
  \else
    \immediate\write\isauxout{\noexpand\expandafter\noexpand\gdef\noexpand\csname #1\noexpand\endcsname{#1}}
    \item{\hbox to \labelwidth{\itemcolor{\hfil#1\ \ }}}
  \fi}

\newcount\chapno
\newcount\chappage
\chapno=0

\newtoks\rightpagehead
\footline={\bgroup{\ifdraft\ifnum\pageno>5\hfill\hbox{\headerfont\today}\fi\fi}\egroup}

\outer\def\section#1{%
  \widow{0.1}%
  \global\subsectno=0%
  \appfalse\secttrue%
  \removelastskip%
  \global\advance\sectno by 1%
  \global\thmno=0\relax%
  \global\randomct=0\relax%
  \edef\numtoks{\ifchap\number\chapno.\fi\number\sectno}%
  \edef\secttitl{#1}%
  \edef\secttitle{Section \numtoks: \secttitl}%
  \rightpagehead={\Black{\headerfont\hfill\secttitle\kern3em\the\pageno}}%
  \vskip 24pt plus 6pt minus 6 pt%
  \ifdraft\let\sectcolor=\OliveGreen\else\let\sectcolor=\Black\fi%
  \message{#1}%
  \ifhypers\hbox{}\fi%
    \futurelet\testchar\maybeoptionsection}

\def\maybeoptionsection{\ifx[\testchar\let\next\optionsection%
	\else\let\next\nooptionsection\fi\next}

\def\optionsection[#1]{%
  \immediate\write\isauxout{\noexpand\forward{\noexpand#1}{Section}{\numtoks}}%
  {\noindent{\draftlabel{#1}\sectionfont\sectcolor{\numtoks}\quad \sectcolor{\secttitl}}}%
  \immediate\write\iscontout{\noexpand\contlist{\noexpand#1}{\secttitl}{\the\pageno}}%
  \htmlanchor{#1}{\makeref{#1}{Section}{\numtoks}}%
  \nobreak\vskip 4ex}

\def\nooptionsection{%
  {\noindent{\sectionfont\sectcolor{\numtoks}\quad \sectcolor{\secttitl}}}%
  \write\iscontout{\noexpand\contlist{0}{\secttitle}{\the\pageno}}%
  \nobreak\vskip 4ex}

\newcount\subsectno
\outer\def\subsection#1{%
  \advance\subsectno by 1%
  \global\randomct=0\relax%
  \vskip 10pt plus 6pt minus 6 pt%
  \widow{.02}%
  \message{#1}%
  \noindent \sectcolor{$\underline{\hbox{\numtoks.\number\subsectno\quad #1}}$}%
  \nobreak\vskip 1ex}


\newif\ifnewgroup\newgroupfalse

%
\def\proclamationsing#1#2#3#4{
  \outer\expandafter\def\csname#1\endcsname{%
  \global\newgroupfalse%
  \ifnum#3<5\global\newgrouptrue\fi%
  \ifnum#3<1\global\newgroupfalse\fi%
  \ifdraft\let\proclaimcolor=\Fuchsia\else\let\proclaimcolor=\Black\fi%
  \gdef\Prnm{#4}%
  \global\advance\thmno by 1%
  \global\randomct=0%
  \ifcase#3
	\stdbreak \or 
	\stdbreak \or 
	\stdbreak \or 
	\stdbreak \or 
	\halfbreak \or 
	\halfbreak \or 
	\halfbreak \or 
	\halfbreak \or 
	\halfbreak \or 
	\halfbreak \or 
	\halfbreak \or 
	\else \stdbreak \fi%
  \edef\proctoks{\ifchap\the\chapno.\fi\ifsect\the\sectno.\fi\the\thmno\ifnum#3=2'\fi\ifnum#3=3''\fi}%
  \widow{0.05}%
  \ifnumct\noindent{\proclaimcolor{%
	\ifnum#3=8*\fi%
	\ifnum#3=9*\fi%
	\ifnum#3=7\hbox to 1ex{$\dag$}\fi%
	\ifnum#3=10\hbox to 0.5ex{$\dag$}\fi%
	\ifnum#3=-1{\bf Important}\ \lowercase\fi%
	\ifnum#3=6\else\ifnum#3=9\else\ifnum#3=10\else{\bf #2}\fi\fi\fi%
	\ifnum#3=6\else\ifnum#3=9\else\ \fi\fi%
	\bf \proctoks}\enspace}%
  \else\noindent{\proclaimcolor{\bf #2}\enspace}%
  \fi%
  \futurelet\testchar\maybeoptionproclaim}}

\def\proclamation#1#2#3{
  \outer\expandafter\def\csname#1\endcsname{%
  \global\newgroupfalse%
  \ifnum#3<5\global\newgrouptrue\fi%
  \ifnum#3<1\global\newgroupfalse\fi%
  \ifdraft\let\proclaimcolor=\Fuchsia\else\let\proclaimcolor=\Black\fi%
  \gdef\Prnm{#2}%
  \global\advance\thmno by 1%
  \global\randomct=0%
  \ifcase#3
	\stdbreak \or 
	\stdbreak \or 
	\stdbreak \or 
	\stdbreak \or 
	\halfbreak \or 
	\halfbreak \or 
	\halfbreak \or 
	\halfbreak \or 
	\halfbreak \or 
	\halfbreak \or 
	\halfbreak \or 
	\else \stdbreak \fi%
  \edef\proctoks{\ifchap\the\chapno.\fi\ifsect\the\sectno.\fi\the\thmno\ifnum#3=2'\fi\ifnum#3=3''\fi}%
  \widow{0.05}%
  \ifnumct\noindent{\proclaimcolor{%
	\ifnum#3=8*\fi%
	\ifnum#3=9*\fi%
	\ifnum#3=7\hbox to 1ex{$\dag$}\fi%
	\ifnum#3=10\hbox to 0.5ex{$\dag$}\fi%
	\ifnum#3=-1{\bf Important}\ \lowercase\fi%
	\ifnum#3=6\else\ifnum#3=9\else\ifnum#3=10\else{\bf #2}\fi\fi\fi%
	\ifnum#3=6\else\ifnum#3=9\else\ \fi\fi%
	\bf \proctoks}\enspace}%
  \else\noindent{\proclaimcolor{\bf #2}\enspace}%
  \fi%
  \futurelet\testchar\maybeoptionproclaim}}

\def\maybeoptionproclaim{\ifx[\testchar\let\next\optionproclaim%
	\else\let\next\nooptionproclaim\fi\next}

\def\optionproclaim[#1]{%
  \bgroup\htmlanchor{#1}{\makeref{\noexpand#1}{\Prnm}{\proctoks}}\egroup%
  \immediate\write\isauxout{\noexpand\forward{\noexpand#1}{\Prnm}{\proctoks}}%
  \draftlabel{#1}%
  \ifnewgroup\bgroup\sl\fi}
\def\nooptionproclaim{\ifnewgroup\bgroup\sl\fi}
\def\endb{\par\stdbreak\egroup\newgroupfalse\randomct=0}

\def\pagelabel#1{%
   \ \unskip%
   \write\isauxout{\noexpand\forward{\noexpand#1}{page}{\the\pageno}}%
   \bgroup\htmlanchor{#1}{\makeref{\noexpand#1}{page}{\the\pageno}}\egroup%
   \draftlabel{page\ #1}%
 }

\def\plot[#1]{
    \ifdraft\let\proclaimcolor=\Fuchsia\else\let\proclaimcolor=\Black\fi%
    \gdef\Prnm{Plot}%
    \global\advance\thmno by 1%
    \ifnumct{\bf \proclaimcolor{Plot\ \proctoks}}%
    \else\noindent{\bf \proclaimcolor{Plot}}%
    \fi%
  \htmlanchor{#1}{\makeref{\noexpand#1}{\Prnm}{\proctoks}}%
  \immediate\write\isauxout{\noexpand\forward{\noexpand#1}{\Prnm}{\proctoks}}%
  \draftlabel{#1}%
  }

\newcount\captionct \captionct=0
\outer\def\caption#1{{\global\advance\captionct by 1%
  \relax%
  \gdef\Irnm{Figure}%
  \edef\numtoks{\ifchap\the\chapno.\fi\ifsect\the\sectno.\fi\the\captionct}%
  \vskip0.5cm\vtop{\noindent{\bf Figure \numtoks.} #1}\vskip0.5cm%
  \futurelet\testchar\maybeoptionproclaim}}
\outer\def\caption#1{{\global\advance\thmno by 1%
  \relax%
  \gdef\Irnm{Figure}%
  \edef\proctoks{\ifchap\the\chapno.\fi\ifsect\the\sectno.\fi\the\thmno}%
  \vskip0.2cm\vtop{\narrower\narrower\noindent{\bf Figure \proctoks.} #1}\vskip0.5cm%
  \futurelet\testchar\maybeoptionproclaim}}

\def\eqlabel#1{
    \ifdraft\let\labelcolor=\Red\else\let\labelcolor=\Black\fi%
    \global\advance\thmno by 1%
    \edef\proctoks{\ifchap\number\chapno.\fi\ifsect\number\sectno.\fi\number\thmno}%
    \outer\expandafter\def\csname#1\endcsname{\proctoks}%
    \htmlanchor{#1}{\makeref{\noexpand#1}{}{\proctoks}}%
    \immediate\write\isauxout{\noexpand\forward{\noexpand#1}{Equation}{(\proctoks)}}%
    \quad\qquad\hfill\eqno\hbox{\hfill\rm\noexpand\labelcolor{(\proctoks)}\draftcmt{\labelcolor{#1}}}%
}
\def\iqlabel#1{
    \ifdraft\let\labelcolor=\Red\else\let\labelcolor=\Black\fi%
    \global\advance\thmno by 1%
    \edef\proctoks{\ifchap\number\chapno.\fi\ifsect\number\sectno.\fi\number\thmno}%
    \outer\expandafter\def\csname#1\endcsname{\proctoks}%
    \htmlanchor{#1}{\makeref{\noexpand#1}{}{\proctoks}}%
    \immediate\write\isauxout{\noexpand\forward{\noexpand#1}{Inequality}{(\proctoks)}}%
    \quad\qquad\hfill\eqno\hbox{\hfill\rm\noexpand\labelcolor{(\proctoks)}\draftcmt{\labelcolor{#1}}}%
}

\def\eqalabel#1{
    \ifdraft\global\let\labelcolor=\Red\else\global\let\labelcolor=\Black\fi%
    \global\advance\thmno by 1%
    \gdef\proctoks{\ifchap\number\chapno.\fi\ifsect\number\sectno.\fi\number\thmno}%
    \outer\expandafter\def\csname#1\endcsname{\proctoks}%
    \htmlanchor{#1}{\makeref{\noexpand#1}{}{\proctoks}}%
    \immediate\write\isauxout{\noexpand\forward{\noexpand#1}{Equation}{(\proctoks)}}%
    \quad\qquad&\hfill\hbox{\hfill\rm\noexpand\labelcolor{(\proctoks)}\draftcmt{\labelcolor{#1}}}%
}

\def\label#1{\optionproclaim[#1]}
\def\label#1{%
    \ifdraft\let\labelcolor=\Red\else\let\labelcolor=\Black\fi%
    \edef\proctoks{\ifchap\the\chapno.\fi\ifsect\the\sectno.\fi\the\thmno}%
    \outer\expandafter\def\csname#1\endcsname{\proctoks}%
    \bgroup\htmlanchor{#1}{\makeref{\noexpand#1}{}{\proctoks}}\egroup%
    \immediate\write\isauxout{\noexpand\forward{\noexpand#1}{Figure}{\proctoks}}%
}

\proclamation{definition}{Definition}{1}
\proclamation{defin}{Definition}{1}
\proclamation{defins}{Definitions}{1}
\proclamation{lemma}{Lemma}{1}
\proclamation{lemmap}{Lemma}{2}
\proclamation{lemmapp}{Lemma}{3}
\proclamation{proposition}{Proposition}{1}
\proclamation{prop}{Proposition}{1}
\proclamation{theorem}{Theorem}{1}
\proclamation{thm}{Theorem}{1}
\proclamation{corollary}{Corollary}{1}
\proclamation{cor}{Corollary}{1}
\proclamation{conjecture}{Conjecture}{1}
\proclamation{proc}{Procedure}{1}
\proclamation{assumption}{Assumption}{1}
\proclamation{axiom}{Axiom}{1}
\proclamation{example}{Example}{0}
\proclamation{examples}{Examples}{0}
\proclamation{examplesnotation}{Examples and notation}{0}
\proclamation{exnoex}{Examples and non-examples}{0}
\proclamation{remark}{Remark}{0}
\proclamation{remarks}{Remarks}{0}
\proclamation{remarksexamples}{Remarks and examples}{0}
\proclamationsing{facts}{Facts}{0}{Fact}
\proclamation{fact}{Fact}{0}
\proclamation{irrelevant}{Irrelevant}{0}
\proclamation{question}{Question}{0}
\proclamation{construction}{Construction}{0}
\proclamation{algorithm}{Algorithm}{0}
\proclamation{problem}{Problem}{0}
\proclamation{table}{Table}{0}
\proclamation{notation}{Notation}{0}
\proclamation{figure}{Figure}{0}
\proclamation{impexample}{Example}{-1}
\proclamation{exercise}{Exercise}{0}
\proclamation{importantexercise}{IMPORTANT EXERCISE}{0}
\proclamation{exercisedag}{Exercise}{7} 
\proclamation{stexercise}{Exercise}{8} 
\proclamation{excise}{Exercise}{6} 
\proclamation{stexcise}{Exercise}{9} 
\proclamation{excisedag}{Exercise}{10} 

\def\doubleindent{\advance \leftskip 2\parindent\advance \rightskip \parindent}
\def\singleindent{\advance \leftskip \parindent\advance \rightskip \parindent}
\def\doubleleftindent{\advance \leftskip 2\parindent}


\def\hpad#1#2#3{\hbox{\kern #1\hbox{#3}\kern #2}}
\def\vpad#1#2#3{\setbox0=\hbox{#3}\vbox{\kern #1\box0\kern #2}}



\def\stack#1#2#3{\vbox{\offinterlineskip
  \setbox2=\hbox{#2}
  \setbox3=\hbox{#3}
  \dimen0=\ifdim\wd2>\wd3\wd2\else\wd3\fi
  \hbox to \dimen0{\hss\box2\hss}
  \kern #1
  \hbox to \dimen0{\hss\box3\hss}}}


\def\hexp#1{%
  \setbox0=\hbox{${}^{#1}$}%
  \hbox to .5\wd0{\box0\hss}}

\def\hsub#1{%
  \setbox0=\hbox{${}_{#1}$}%
  \hbox to .5\wd0{\box0\hss}}


\def\mOth{\mathsurround=0pt}
\newdimen\pOrenwd \setbox0=\hbox{\tenex B} \pOrenwd=\wd0
\def\leftextramatrix#1{\begingroup \mOth
  \setbox0=\vbox{\def\cr{\crcr\noalign{\kern2pt\global\let\cr=\endline}}
    \ialign{$##$\hfil\kern2pt\kern\pOrenwd&\thinspace\hfil$##$\hfil
      &&\quad\hfil$##$\hfil\crcr
      \omit\strut\hfil\crcr\noalign{\kern-\baselineskip}
      #1\crcr\omit\strut\cr}}
  \setbox2=\vbox{\unvcopy0 \global\setbox1=\lastbox}
  \setbox2=\hbox{\unhbox1 \unskip \global\setbox1=\lastbox}
  \setbox2=\hbox{$\kern\wd1\kern-\pOrenwd \left[ \kern-\wd1
    \vcenter{\unvbox0 \kern-\baselineskip} \,\right]$}
  \null\;\vbox{\box2}\endgroup}


\def\COMMENT#1\par{\bigskip\hrule\smallskip#1\smallskip\hrule\bigskip}

{\catcode`\^^M=12 \endlinechar=-1 %
 \gdef\xcomment#1^^M{\def\test{#1}
   \ifx\test\endcomment \let\next=\endgroup
   \else \let\next=\xcomment \fi
   \next}
}
\def\dospecials{\do\ \do\\\do\{\do\}\do\$\do\&%
  \do\#\do\^\do\^^K\do\_\do\^^A\do\%\do\~}
\def\makeinnocent#1{\catcode`#1=12}
\def\comment{\begingroup
  \let\do=\makeinnocent \dospecials
  \endlinechar`\^^M \catcode`\^^M=12 \xcomment}
{\escapechar=-1
 \xdef\endcomment{\string\\endcomment}
}



\def\draftcmt#1{\ifdraft\llap{\smash{\raise2ex\hbox{{#1}}}}\fi}
\def\draftind#1{\llap{\smash{\raise-1.3ex\hbox{{#1}}}}}
\def\draftlabel#1{\ifdraft\draftcmt{\footfonts\Green{{\string#1}\ }\unskip\fi}}
\def\draftlabel#1{%
  \ifx#1\empty\else
  \ifdraft\draftcmt{\footttt\Green{\string #1}\ \ }\unskip\fi\fi}

\newbox\strikebox
\def\strike#1{\setbox\strikebox=\hbox{#1}%
  \rlap{\raise0.4ex\hbox to \wd\strikebox{\leaders\hrule height 0.8pt depth 0pt\hfill}}%
  #1}
\def\strikem#1{\setbox\strikebox=\hbox{$#1$}%
  \rlap{\raise0.4ex\hbox to \wd\strikebox{\leaders\hrule height 0.8pt depth 0pt\hfill}}%
  \hbox{$#1$}}

\newcount\hours
\newcount\minutes
\hours\time \divide\hours 60
\minutes-\hours \multiply\minutes 60\advance\minutes\time
\edef\timehhmm{\ifnum\hours<10 0\fi\the\hours
:\ifnum\minutes<10 0\fi\the\minutes}

\footline={{\rm\hfill\folio\hfill}}

\def\cite{\futurelet\testchar\maybeoptioncite}
\def\maybeoptioncite{\ifx[\testchar \let\next\optioncite
	\else \let\next\nooptioncite \fi \next}

\def\nooptioncite#1{\expandafter\ifx\csname#1\endcsname\relax
  {\nodefmsg{#1}\bf`#1'}\else
  \htmllocref{#1}{[\csname#1\endcsname]}{}\fi}

\def\optioncite[#1]#2{\ifx\csname#1\endcsname\relax\else
  \htmllocref{#2}{[\csname#2\endcsname, #1]}{}\fi}


\def\strutdepth{\dp\strutbox}
\def\strutdepth{4.3ex}
\def\margin#1{{\strut\vadjust{\kern-\strutdepth%
        \vtop to\strutdepth{\baselineskip\strutdepth%
        \llap{\vbox{\hsize=3em\noindent #1}\hskip1em}\null}}}}
\def\rmargin#1{\strut\vadjust{\kern-\strutdepth%
        \vtop to\strutdepth{\baselineskip\strutdepth%
        \vss\rlap{\ \hskip\hsize\ \vtop{\hsize=4em\noindent #1}}\null}}}
\def\warning{\vskip0ex\noindent\llap{\lower4pt\hbox to 1.0cm{\psline(0,.3)(.27,.3)(0,.09)(.03,0)(.27,0)(.3,.03)(.27,.06)}\hfill\ }\indent}
\def\warning{\margin{\lower4pt\hbox to 0.6cm{\psline(0,.3)(.27,.3)(0,.09)(.03,0)(.27,0)(.3,.03)(.27,.06)}\hfill\ }}
\def\warning{\margin{\hbox to 0.6cm{\psline(0,.3)(.27,.3)(0,.09)(.03,0)(.27,0)(.3,.03)(.27,.06)}\hfill\ }}
\def\warning{\noindent{\bf Warning: }\margin{\ \ \ \ \lower-4pt\hbox to 0.6cm{\psline[linewidth=2pt](0,0)(.5,0)(0.53,-0.02)(0.55,-0.05)(0.53,-0.1)(0.03,-0.35)(0.01,-0.4)(0.03,-0.43)(0.06,-0.45)(0.56,-0.45)}\hfill}}
\def\warningdown#1{\noindent{\bf Warning: }\bgroup\def\strutdepth{#1}\margin{\ \ \
\ \hbox to 0.6cm{\psline[linewidth=2pt](0,0)(.5,0)(0.53,-0.02)(0.55,-0.05)(0.53,-0.1)(0.03,-0.35)(0.01,-0.4)(0.03,-0.43)(0.06,-0.45)(0.56,-0.45)}\hfill}\egroup}
\def\quote{\bgroup\footfont\parindent=0pt%
	\baselineskip=8pt\hangindent=15em\hangafter=0}
\def\endquote{\egroup}

\def\bmargin#1{{\strut\vadjust{\kern-\strutdepth%
        \vtop to\strutdepth{\baselineskip\strutdepth%
        \llap{\vbox{\hsize=3em\noindent #1}\hskip2em}\null}}}}

\newdimen\hngind \hngind=\parindent \advance\hngind by 3.9em\relax
\def\nonuritem{{\vskip0pt\hangafter=0\global\hangindent\hngind%
        {\vphantom{(\the\randomct)\hskip.5em}}}}
\def\nonuritem{{\parindent=0pt\hngind=\parindent \advance\hngind by 1.9em%
        \vskip0pt\hangafter=0\global\hangindent\hngind%
        \vphantom{{(\the\randomct)\hskip.5em}}}}
\def\nonuritem{{\idoitem{}}}
\def\pitem{{\global\advance\randomct by 1%
	\vskip1pt\hangafter=0\global\hangindent\hngind%
	{{\the\randomct.\hskip.5em}}}}

%
\def\idoitem#1{{\parindent=0pt\hngind=\parindent \advance\hngind by 3.7em%
        \vskip0pt\hangafter=0\global\hangindent\hngind%
        {$\hbox to3.7em{\hfill#1\hskip.5em}$}}}%
\def\romanitem{\global\advance\randomct by 1%
    \edef\itoks{(\romannumeral\randomct)}%
    \edef\Irnm{}%
    \idoitem{\itoks}%
    \futurelet\testchar\maybeoptioniitem} 
\def\sitem{\global\advance\randomct by 1%
    \edef\itoks{(\the\randomct)}%
    \edef\Irnm{Step}%
    \idoitem{\itoks}%
    \futurelet\testchar\maybeoptioniitem} 
\def\Rsitem{\global\advance\rrandomct by 1%
    \edef\itoks{(\uppercase\expandafter{\romannumeral\rrandomct)}}%
    \edef\Irnm{Step}%
    \idoitem{\itoks}%
    \futurelet\testchar\maybeoptioniitem} 
\def\citem{\global\advance\randomct by 1%
    \edef\itoks{(\the\randomct)}%
    \edef\Irnm{Condition}%
    \idoitem{\itoks}%
    \futurelet\testchar\maybeoptionitem} 
\def\iitem{\global\advance\randomct by 1%
    \edef\itoks{(\the\randomct)}%
    \edef\Irnm{}%
    \idoitem{\itoks}%
    \futurelet\testchar\maybeoptioniitem} 
\def\rpitem{\global\advance\randomct by 1%
    \edef\itoks{(\the\randomct')}%
    \edef\Irnm{}%
    \idoitem{\itoks}%
    \futurelet\testchar\maybeoptioniitem} 

\def\maybeoptioniitem{\ifx[\testchar\let\next\optioniitem%
	\else\let\next\nooptioniitem\fi\next}
\def\optioniitem[#1]{%
  \bgroup\htmlanchor{#1}{\makeref{\noexpand#1}{\Irnm}{\itoks}}\egroup%
  \immediate\write\isauxout{\noexpand\forward{\noexpand#1}{\Irnm}{\itoks}}%
  \draftlabel{#1}}
\def\nooptioniitem{}
\def\bulletpoint{{\parindent=0pt\hngind=\parindent \advance\hngind by 1.2em%
        \vskip0pt\hangafter=0\global\hangindent\hngind%
        {$\hbox to1.2em{\hfill$\bullet$\hskip.5em}$}}}%
\def\bulletpoint{{\hngind=\parindent \advance\parindent by -1.2em%
        \vskip0pt\hangafter=0\global\hangindent\hngind%
        {$\hbox to1.2em{\hfill$\bullet$\hskip.5em}$}}}%
\def\nopoint{{\parindent=0pt\hngind=1cm \advance\hngind by 1.2em%
        \vskip0pt\hangafter=0\global\hangindent\hngind%
        {$\hbox to1.2em{\hfill\ \hskip.5em}$}}}%
\def\nopoint{{\hngind=\parindent \advance\parindent by -1.2em%
        \vskip0pt\hangafter=0\global\hangindent\hngind%
        {$\hbox to1.2em{\hfill\ \hskip.5em}$}}}%
\def\bbulletpoint{{\bgroup\parskip=0.1cm\parindent=1.5cm\hngind=\parindent%
	\advance\hngind by 1.5em%
        \vskip0pt\hangafter=0\global\hangindent\hngind%
        {$\hbox to1.5em{\hfill$\bullet$\hskip.5em}$}\egroup}}%
\def\bnobulletpoint{{\bgroup\parskip=0.1cm\parindent=1.5cm\hngind=\parindent%
	\advance\hngind by 1.5em%
        \vskip0pt\hangafter=0\global\hangindent\hngind%
        {$\hbox to1.5em{\hfill\hskip.5em}$}\egroup}}%

\let\ritem=\iitem

\def\rlabel#1{{%
  \bgroup\htmlanchor{#1}{\makeref{\noexpand#1}{item}{\the\randomct}}\egroup%
  \immediate\write\isauxout{\noexpand\forward{\noexpand#1}{item}{\the\randomct}}%
  \draftlabel{item\ #1}}}

\def\idoitemitem#1{{\hngind=\parindent \advance\hngind by 4.3em%
        \vskip-\parskip\hangafter=0\global\hangindent\hngind%
        {$\hbox to3.9em{\hfill#1\ }$}}}
\def\idoitemitemm#1{{\hngind=\parindent \advance\hngind by 7.3em%
        \vskip-\parskip\hangafter=0\global\hangindent\hngind%
        {$\hbox to6.9em{\hfill#1\ }$}}}
\let\riitem=\idoitemitem

\def\thmparens#1{\ignorespaces{\rm(#1)}}

\newdimen\thmskip \thmskip=0.7ex


\newdimen\exerciseindent \exerciseindent=4.5em
\newdimen\exerciseitemindent \exerciseitemindent=1.9em
\newdimen\exitindwk \exitindwk=3.7em
\def\exitem{\global\advance\randomct by 1%
	\vskip0pt\indent\hangindent\exitindwk \hskip\exerciseitemindent%
	\llap{\romannumeral\randomct)\hskip0.5em}\ \hskip-0.5em}
\def\stexitem{\global\advance\randomct by 1%
	\par\indent\hangindent\exitindwk \hskip\exerciseitemindent%
	\llap{(\romannumeral\randomct)\hskip-.1ex*\enspace}\ignorespaces}


\def\exercises{\global\thmno=0\vskip0.9cm%
	\noindent{\bf Exercises for Section~\ifchap\the\chapno.\fi\the\sectno}}
\def\solutions{\global\thmno=0\bigskip\widow{0.05}%
	\noindent{\bf Solutions for Section~\ifchap\the\chapno.\fi\the\sectno}}
\def\solution#1{\removelastskip\vskip 5pt%
	\global\randomct=0%
	\parindent=\exerciseindent%
	\par\hangindent\exerciseindent\indent\llap{\hbox to \exerciseindent{%
	  {\bf \unskip\ref{#1}}}}}
\def\solution#1{\removelastskip\vskip 5pt\par\noindent {\bf \refn{#1}:}}


\newdimen\wddim
\def\edinsert#1{\setbox0=\hbox{#1}\wddim=\wd0\divide\wddim by 2%
\smash{\hskip-\wddim\raise15pt\hbox to 0pt{\Red{$\underbrace{\box0}$}}\hskip\wddim}}

\def\midline#1{\raise 0.20em\hbox to #1{\vrule depth0pt height 0.405pt width #1}}

\def\buildrelu#1\over#2{\mathrel{\mathop{\kern0pt #1}\limits_{#2}}}



\def\bmatrix#1#2{{\left(\vcenter{\halign
  {&\kern#1\hfil$##\mathstrut$\kern#1\cr#2}}\right)}}

\def\rightarrowmat#1#2#3{
  \setbox1=\hbox{\small\kern#2$\bmatrix{#1}{#3}$\kern#2}
  \,\vbox{\offinterlineskip\hbox to\wd1{\hfil\copy1\hfil}
    \kern 3pt\hbox to\wd1{\rightarrowfill}}\,}

\def\leftarrowmat#1#2#3{
  \setbox1=\hbox{\small\kern#2$\bmatrix{#1}{#3}$\kern#2}
  \,\vbox{\offinterlineskip\hbox to\wd1{\hfil\copy1\hfil}
    \kern 3pt\hbox to\wd1{\leftarrowfill}}\,}

\def\rightarrowbox#1#2{
  \setbox1=\hbox{\kern#1\hbox{\small #2}\kern#1}
  \,\vbox{\offinterlineskip\hbox to\wd1{\hfil\copy1\hfil}
    \kern 3pt\hbox to\wd1{\rightarrowfill}}\,}

\def\leftarrowbox#1#2{
  \setbox1=\hbox{\kern#1\hbox{\small #2}\kern#1}
  \,\vbox{\offinterlineskip\hbox to\wd1{\hfil\copy1\hfil}
    \kern 3pt\hbox to\wd1{\leftarrowfill}}\,}


\newcount\countA
\newcount\countB
\newcount\countC

\def\monthname{\begingroup
  \ifcase\number\month
    \or January\or February\or March\or April\or May\or June\or
    July\or August\or September\or October\or November\or December\fi
\endgroup}

\def\dayname{\begingroup
  \countA=\number\day
  \countB=\number\year
  \advance\countA by 0 
  \advance\countA by \ifcase\month\or
    0\or 31\or 59\or 90\or 120\or 151\or
    181\or 212\or 243\or 273\or 304\or 334\fi
  \advance\countB by -1995
  \multiply\countB by 365
  \advance\countA by \countB
  \countB=\countA
  \divide\countB by 7
  \multiply\countB by 7
  \advance\countA by -\countB
  \advance\countA by 1
  \ifcase\countA\or Sunday\or Monday\or Tuesday\or Wednesday\or
    Thursday\or Friday\or Saturday\fi
\endgroup}

\def\timename{\begingroup
   \countA = \time
   \divide\countA by 60
   \countB = \countA
   \countC = \time
   \multiply\countA by 60
   \advance\countC by -\countA
   \ifnum\countC<10\toks1={0}\else\toks1={}\fi
   \ifnum\countB<12 \toks0={\sevenrm AM}
     \else\toks0={\sevenrm PM}\advance\countB by -12\fi
   \relax\ifnum\countB=0\countB=12\fi
   \hbox{\the\countB:\the\toks1 \the\countC \thinspace \the\toks0}
\endgroup}

\def\timestamp{\dayname, \the\day\ \monthname\ \the\year, \timename}

%

%


\def\frac#1#2{{#1 \over #2}}

\let\text\hbox
\def\overset#1\to#2{\ \mathop{\buildrel #1 \over #2}\ }
\def\underset#1\to#2{{\mathop{\buildrel #1 \over #2}}}

\def\myoperator{NOOP}
\def\mymathopsp{\futurelet\testchar\maybesubscriptop}
\def\maybesubscriptop{\ifx_\testchar \let\next\subscriptop
	\else \let\next\nosubscriptop \fi \next}
\def\subscriptop_#1{{\rm \myoperator}_{#1}}
\def\nosubscriptop{\futurelet\testchar\maybeparenop}
\def\maybeparenop{{\mathop{\rm \myoperator}\nolimits%
	\ifx(\testchar \relax\else\,\fi}}
\font\tenmsb=msbm10
\font\sevenmsb=msbm10 at 7pt
\font\fivemsb=msbm10 at 5pt
\newfam\msbfam
\textfont\msbfam=\tenmsb
\scriptfont\msbfam=\sevenmsb
\scriptscriptfont\msbfam=\fivemsb

\def\hexnumber#1{\ifcase#1 0\or1\or2\or3\or4\or5\or6\or7\or8\or9\or
	A\or B\or C\or D\or E\or F\fi}

\mathchardef\subsetneq="2\hexnumber\msbfam28



\def\dim{\gdef\myoperator{dim}\mymathopsp}

\def\notdiv{\hbox{$\not\hbox to -0.2ex{$|$}$}\hskip0.5em}

\def\Ext{\mathop{\rm Ext}\nolimits}

\def\Ass{\mathop{\rm Ass}\nolimits}

\def\projdim{\mathop{\rm pd}\nolimits}

\def\min{\mathop{\rm min}}
\def\max{\mathop{\rm max}}

\def\depth{\mathop{\rm depth}\nolimits\,}

\def\mod{\mathop{\rm mod}\nolimits\,}


\def\deg{\mathop{\rm deg}\nolimits}
\def\depth{\mathop{\rm depth}\nolimits}

\def\dim{\mathop{\rm dim}\nolimits}

\def\Hom{\mathop{\rm Hom}\nolimits}

\def\Ext{\mathop{\rm Ext}\nolimits}

\def\semidirect{{\mathop{\hbox to 1.1em{\hskip.5em\psline[linewidth=0.3pt](-.18,+.01)(.18,.15)(.18,.01)(-.18,.15)}}}}


\def\spm{\raise0.5pt\hbox{\bgroup\textfont2=\fivesy\relax$\pm$\egroup}}
\def\dotdiv{\mathop{\hbox{\rlap{\raise-0.1em\hbox{--}}\hskip0.1em\raise0.7ex\hbox{.}}\ }}

\def\uspace{\hbox{$\underline{\hbox to 1ex{\ \hfil}}$}}

\def\lineunderblank{\underline{\ \hbox to 0.5ex{$\vphantom{a}$\hfill\ }}}
\def\lineunderblankb#1{\underline{\ \hbox to #1{$\vphantom{a}$\hfill\ }}}
\def\lineaboveblank{\overline{\ \hbox to 0.5ex{$\vphantom{a}$\hfill\ }}}
\def\lineaboveblank{\overline{\ \hbox to 0.5ex{$\vphantom{a}$\hfill\ }}}
\def\blank#1{\hbox{$\underline{\hbox to #1em{\ \hfill}}$}}


\def\qedbox{\hbox{\vbox{\hrule\hbox{\vrule\kern3pt\vbox{\kern6pt}\kern3pt\vrule}\hrule}}}
\def\qed{\ \unskip\hfill\qedbox\vskip\thmskip} 


\def\today{\ifcase\month\or January\or February\or March\or
  April\or May\or June\or July\or August\or September\or
  October\or November\or December\fi
  \space\number\day, \number\year}\relax


\newwrite\isauxout
\openin1\jobname.aux
\ifeof1\message{No file \jobname.aux}
       \else\closein1\relax\input\jobname.aux
       \fi
\immediate\openout\isauxout=\jobname.aux
\immediate\write\isauxout{\relax}

\newwrite\iscontout
\immediate\openout\iscontout=\jobname.cont
\immediate\write\iscontout{\relax}

\def\contlist#1#2#3{\ifchap\vskip1pt\noindent \ \ \ifx0#1#2\hfill#3%
		\else\ref{#1}: #2\hfill#3%
                \fi%
	\else \vskip1pt\noindent\ref{#1}: #2\hfill#3\fi}
\def\contlist#1#2#3{\vskip1pt\noindent \ \hskip1em%
	\ifchap\ifx0#1#2\hfill#3%
		\else\ref{#1}: #2\hfill#3%
                \fi%
	\else \ref{#1}: #2\hfill#3\fi}
\def\contnosectlist#1#2#3{\vskip1pt\noindent\hskip1em#2\hfill#3}


%
\newwrite\isindexout
\immediate\openout\isindexout=\jobname.index
\font\footttt=cmtt6

\def\padno{
\ifnum\pageno<10 00\the\pageno\else\ifnum\pageno<100 0\the\pageno\else\the\pageno\fi\fi}

\def\index#1{\rlap{\ifdraft\draftind{{\footttt\Blue{#1\ \ }}}\unskip\fi%
\write\isindexout{\noexpand\indexentry {#1}{\padno}}}\unskip\ignorespaces}

\def\indexn#1{\rlap{\ifdraft\draftind{{\footttt\Brown{#1\ \ }}}\unskip\fi%
\write\isindexout{\noexpand\indexentryn {#1}{\padno}}}\unskip\ignorespaces}

\def\obeyspaces{\catcode`\ =\active\catcode`\	=\active}
{\obeyspaces\global\let =\space} 
{\obeyspaces\gdef {\hskip.5em}\gdef	{\hskip3em}}
{\catcode`\^^M=\active %
\gdef\obeylines{\catcode`\^^M=\active \gdef^^M{\vskip.0ex\ }}}%

\font\stt=cmtt10
\def\dMaccode{\vskip3ex\bgroup\advance\leftskip1ex\parindent=0pt%
	\baselineskip=9pt\obeylines\obeyspaces%
	\catcode`\^=11\catcode`\_=11\catcode`\#=11\catcode`\%=11%
	\catcode`\{=11\catcode`\}=11%
	\stt}
\def\enddMaccode{\removelastskip\unskip\egroup}
\def\Maccode{\bgroup \catcode`\^=11 \catcode`\_=11 \stt"}
\def\endMaccode{\unskip"\egroup}
\def\Maccode{\bgroup \catcode`\^=11 \catcode`\_=11 \stt}
\def\endMaccode{\unskip\egroup}

\gdef\poetry{\vskip0pt\bgroup\parindent=0pt\obeylines\obeyspaces}
\gdef\poetrycenter{\vskip0pt\bgroup\parskip=0pt\leftskip=0pt plus 1fil\rightskip=0pt plus
1fil\parfillskip=0pt\obeylines}
\gdef\endpoetry{\vskip0pt\egroup}

\magnification=1200

\newwrite\isauxout
\openin1\jobname.aux
\ifeof1\message{No file \jobname.aux}
       \else\closein1\relax\input\jobname.aux
       \fi
\immediate\openout\isauxout=\jobname.aux
\immediate\write\isauxout{\relax}

\def\BaHH{{\mathop{\rm BHH}\nolimits}}

\abovedisplayskip=2pt%
\belowdisplayskip=2pt%

\sectfalse
\drafttrue
\draftfalse
\colortrue
\ifcolor\input colordvi\else\input blackdvi\fi%
\ifdraft\footline={{\rm \hfill page \folio -- \today\ --
\timehhmm\hfill}}\else
\ifnum\folio>1\footline={{\rm\hfill\folio\hfill}}\else\footline={{\rm\hfill\folio\hfill}}\fi\fi

\centerline{\sectionfont
Fluctuations in depth and associated primes of powers of ideals}
\vskip 2mm
\centerline{\sectionfont Roswitha Rissner and Irena Swanson}
\centerline{\tt Roswitha.Rissner@aau.at, irena@purdue.edu}

\medskip
\bgroup
\narrower\narrower
\noindent
{\bf Abstract}
We count the numbers of associated primes of powers of ideals
as defined in~\cite{BHH}.
We generalize those ideals to monomial ideals
$\BaHH(m,r,s)$ for $r \ge 2$, $m, s \ge 1$;
we establish partially the associated primes of powers of these ideals,
and we establish completely
the depth function of quotients by powers of these ideals:
the depth function is periodic of period~$r$ repeated $m$ times on the initial interval
before settling to a constant value.
The number of needed variables for these depth functions
are lower than those from general constructions in~\cite{HNTT}.

\egroup
\medskip

This paper was motivated by results from
Herzog and Hibi~\cite{HH}
and Bandari, Herzog and Hibi \cite{BHH}
that construct monomial ideals $I$
with various properties of the depth function $n \mapsto \depth(R/I^n)$.
In particular,
Herzog and Hibi \cite{HH}
constructed for any non-increasing eventually constant sequence $\{a_n\}$
a monomial ideal such that for all integers~$n$,
$\depth(R/I^n) = a_n$.
In general,
the depth function need not be monotone,
as shown by an example in~\cite{HH}.
Bandari, Herzog and Hibi \cite{BHH}
constructed for each positive integer~$m$ a monomial ideal $I$
for which the depth function
takes on values $0, 1$, repeated $m$ times,
followed by $0$ and then by constant~$2$.
Thus this function has a global maximum,
exactly $m$ strict local maxima
and exactly $m+1$ strict local minima.
This was the first example of prescribed depth periodicity of period~$2$
on a segment of the domain.
We point out that
a later paper, \cite{HNTT}, by H\`a, Nguyen, Trung and Trung,
establishes more generally
for any eventually constant $\hbox{\tenbbold N}_0$-valued sequence~$\{a_n\}$
the existence of a monomial ideal~$Q$ in a polynomial ring~$S$
satisfying $\depth(S/Q^n) = a_n$ for all~$n$.
This completely determines all depth functions of powers of ideals.

Part of our long-term goal is to shed light similarly
on the possible functions $n \mapsto \#\Ass(R/I^n))$
for ideals~$I$ in Noetherian rings~$R$.
Certainly these functions are all positive-integer valued
and eventually constant by a result of Brodmann~\cite{Brodmann}.
The second author and Weinstein proved in~\cite{WS}
that for every non-increasing sequence $\{a_n\}$ of positive integers
there exists a family of monomial ideals $I$
such that for all~$n$,
the number of associated primes of $I^n$ is $a_n$.
For arbitrary (necessarily eventually constant) sequences of positive integers
much less is known.
If some $a_n$ equals~$1$,
then if we are to vary over monomial ideals
it is necessary that all $a_m$ for $m \ge n$ also be equal to~$1$.
If we do not restrict to monomial ideals,
then a big jump can occur from $a_1$ to $a_2$ even if we restrict to prime
ideals;
a result from~\cite{KS} proves that $a_2$ is not bounded above
by any polynomial function in the number of variables in the ring.

We present in \ref{thmBHH2count}
the function $n \mapsto \#\Ass(R/I^n))$ for ideals~$I$ introduced
by Bandari, Herzog and Hibi in~\cite{BHH}.
Once we completed the count of all associated primes
and observed certain partial periodicity of period~$2$,
we introduced a more general family of ideals, $\BaHH(m,r,s)$ with $r \ge 2$,
$m, s \ge 1$;
in this notation,
the original Bandari--Herzog--Hibi ideals are $\BaHH(m,2,2)$,
and in \ref{thmBHH2count}
we count more generally the associated primes of $\BaHH(m,2,s)^n$ as:
$$
(3-\delta_{1=n})^m +
\left( \sum_{\ell=0}^m \sum_{t=b(\ell)}^m { {m}\choose{\ell}} {{\ell}\choose{\ell+t-m} } \right) +
\cases{0, & if $n \le 2m$ and $n$ is even; \cr
1, & otherwise. \cr }
$$
where $b(\ell) = \max\{n-1-\ell,m-\ell\}$
and $\delta_C$ equals $1$ if the condition~$C$ is true and $0$ otherwise.
In particular,
the number of associated primes of $\BaHH(m,2,s)$ is $2^m + 3^m + 1$,
the number of associated primes of $\BaHH(m,2,s)^2$ is $2 \cdot 3^m$,
and when $n \ge 2m+2$,
the number of associated primes is $1 + 3^m$.
In \ref{thmcountmax} we prove that
the function $n \mapsto \#\Ass(R/\BaHH(m,2,s)^n)$
has exactly $\left\lceil {m-1 \over 2} \right\rceil$ local maxima.
The global maximum $2 \cdot 3^m + 1$
is achieved exactly at
$n = 3, 5, \ldots, 2 \left\lceil {m-1 \over 2} \right\rceil + 1$.
This function is periodic of period~2
when restricted to $[3, 2 \left\lceil {m-1 \over 2} \right\rceil + 1]$.

We present this count of associated primes in two different ways.
As a result,
\ref{rmknewcomb} proves an identity of binomial expressions
that we have not found in the literature.

For $r > 2$ we completely describe and count all the associated primes
that contain one of the special variables $c_1, \ldots, c_s$ (and hence all),
and we give some properties and descriptions of the associated primes
that do not contain these special variables.
The latter associated primes satisfy persistence,
namely that if a prime ideal not containing the special variable
is associated to an $n$th power,
then it is associated to all higher powers.
This persistence is not in general satisfied by the associated primes containing
the special variable.

Seidenberg proved in \cite[Point 65]{Sei74}
that there exists a primitive recursive function $B(n, d)$
such that any ideal $I$ in a polynomial ring
in $n$ variables over a field with generators of degree at most $d$
has at most $B(n,d)$ associated primes.
Ananyan and Hochster \cite{AH} proved
that there exists a primitive recursive function $E(g,d)$
such that any ideal $I$ in a polynomial ring over a field
with at most $g$ generators of degrees $d$ or less
has at most $E(g,d)$ associated primes.
The ideal $\BaHH(m,2,1)$
is in a polynomial ring with $2m+3$ variables
and has $2m + 5$ generators of degrees up to~$9$.
Its third power has $2 \cdot 3^m + 1$ associated primes,
and by \ref{lmpoweregen} it has at most
${2m+1 \choose 2} + {2m + 6 \choose 3}$ generators
of degrees up to~$23$.
For large~$m$, this number of generators is less than or equal to $2m^3$,
showing that for large $m$ and $n$,
$$
\eqalignno{
&B(2m+3,23) \ge 2 \cdot 3^m + 1,
\hbox{ i.e., }
B(n,23) \ge 2 \cdot (\sqrt 3)^{n - 3} + 1, \cr
&E(2m^3,23) \ge 2 \cdot 3^m + 1,
\hbox{ i.e., }
E(n, 23) \ge
2 \cdot 3^{\root 3 \of {n/2}} + 1.
}
$$
Asymptotically,
the lower bound here for $B(n,d)$ is stronger
than the bound $3^{n/3}$ in \cite{KS},
but the lower bound for $E(n,d)$ here is weaker
than the bound $3^{\sqrt{2n}-1}$ in \cite{KS}.

In \ref{thmdepth2}
we prove that the function $n \mapsto \depth(R/\BaHH(m,r,s)^n)$
is periodic of period $r$ when restricted to the interval $[1, \ldots, rm+1]$,
that it has exactly $m+1$ local minima, all on that interval and equal to~$0$,
that all other values on that interval are~1,
and that the only value outside of the interval is~$s$.
More generally,
we show that
an $e$-fold splitting of $\BaHH(m,r,s)$ gives 
an ideal $I$ in a ring $A$ such that
$$
\depth\left({A \over I^n}\right)
= \cases{
e-1, & if $n = ru + 1$ with $u = 0, \ldots, m$; \cr
e, & if $n \le rm+1$ and $n \not \equiv 1 \mod r$; \cr
s+e-1, & otherwise, i.e., if $n > mr+1$. \cr
}
$$
We point out that
the construction in~\cite{HNTT} by H\`a, Nguyen, Trung and Trung
of the monomial ideal with the same depth function uses at least
$e + 4 (rm-m) + 3 s$ variables,
whereas our construction uses $rm + r + s + e - 1$.
The difference
$3rm - 4m - r + 2s + 1 = 2(r-2)m + (m-1)r + 2s + 1$
in the number of variables
is always positive since $r \ge 2$ (for periodicity) and $m, s \ge 1$.

\vskip 3mm
{\bf Acknowledgement.}
We thank Ram Goel for his contributions
to the understanding of associated primes,
including for writing a Python program to compute our projected numbers.
His work was crucial to helping us correct our initial patterns.

\section{Generalized Bandari--Herzog--Hibi ideals}

\defin[definI]%
Let $m, r$ and $s$ be positive integers with $r \ge 2$.
Let $c_1, \ldots, c_s, a_j, x_{i,j}$ be variables over a field $k$
where $i \in [m]$ and $j \in [r]$.
For the sake of notation,
we identify $a_{r+1} = a_1$, $a_0 = a_r$,
$x_{i,r+1} = x_{i,1}$,
$x_{i,0} = x_{i,r}$ for all $i$,
and more generally, $x_{i,j} = x_{i,(j \mod r)}$.
We define
$$
\eqalignno{
B_0(r) &= (a_j^6, a_j^5 a_{j+1} : j = 1, \ldots, r), \cr
B_c(r,s) &= (c_1, \ldots, c_s) a_1^4 a_2^4 \cdots a_r^4; \cr
X(m,r) &= (a_j^4 x_{i,j} x_{i,j+1}^2 : i \in [m], j \in [r]), \cr
\BaHH(m,r,s) &= B_0(r) + B_c(r,s) + X(m,r). \cr
}
$$
We call $\BaHH(m,r,s)$ the {\bf Bandari-Herzog-Hibi} ideals.
When $s = 1$,
we write $c = c_1$.

When $m, r$ and $s$ are clear,
we write these ideals as $B_0, B_c, X, B$, respectively.

We will refer to elements $a_j^4 x_{i,j} x_{i,j+1}^2$ as $h_{i,j}$.
\endb

We name these ideals in honor of the Bandari-Herzog-Hibi paper~\cite{BHH}
which originated the ideals $\BaHH(m,2,2)$.

Understanding the associated primes of powers of these ideals
is important for understanding the depth function of their quotients.
For the primary decompositions part,
we prove in \ref{thmspreadassoc}
that it suffices to find the decompositions in case $s = 1$.
This reduction greatly simplifies the notation
and speeds up any concrete calculations of the associated primes
and thus the counting.

\lemma[lmpoweregen]
Let $J_1, J_2, J_3$ be ideals in a ring $R$ such that 
$(J_1 + J_2)^2 \subseteq J_2^2$.
Then for all positive integers $n$,
$(J_1 + J_2 + J_3)^n = J_1 J_3^{n-1} + (J_2 + J_3)^n$.

Thus with $B = \BaHH(m,r,s)$,
$B^n = B_c X^{n-1} + (B_0 + X)^n = B_c B^{n-1} + (B_0 + X)^n$.
\endb

\proof
The first display holds trivially for $n = 1$.
The equality
$(J_1 + J_2 + J_3)^2
= J_1^2  + J_1 J_2 + J_1 J_3 + J_2^2 + J_2 J_3 + J_3^2$
$= J_1 J_3 + J_2^2 + J_2 J_3 + J_3^2$
$= J_1 J_3 + (J_2 + J_3)^2$
proves the case $n = 2$.
Then by induction on $n \ge 2$,
$$
\eqalignno{
(J_1 + J_2 + J_3)^n
&= (J_1 + J_2 + J_3) (J_1 + J_2 + J_3)^{n-1} \cr
&= (J_1 + J_2 + J_3) \left(J_1 J_3^{n-2} + (J_2 + J_3)^{n-1}\right) \cr
&= J_1^2 J_3^{n-2} + J_1 (J_2 + J_3)^{n-1} + J_1 J_2 J_3^{n-2}
+ (J_2 + J_3)^n + J_1 J_3^{n-1}
\cr
&= J_1 (J_2 + J_3)^{n-1} + (J_2 + J_3)^n \cr
&= J_1 \sum_{i=0}^{n-1} J_2^i J_3^{n-1-i} + (J_2 + J_3)^n \cr
&= J_1 J_3^{n-1} + \sum_{i=1}^{n-1} J_1 J_2^i J_3^{n-1-i} + (J_2 + J_3)^n \cr
&\subseteq J_1 J_3^{n-1} + (J_2 + J_3)^n. \cr
}
$$
Since also the last ideal is contained in the first in this display,
the conclusion follows.

The second part follows with $J_1 = B_c$, $J_2 = B_0$ and $J_3 = X$,
since $J_1^2 \subseteq (a_1^6) (a_2^6) \subseteq J_2^2$
and $J_1 J_2 \subseteq (a_j^6, a_j^5 a_{j+1}: j \in [r]) (a_1^4 \cdots a_r^4)
\subseteq ((a_j^5 a_{j+1})^2, (a_j^6) (a_{j+1}^5 a_{j+2}): j \in [r])
\subseteq J_2^2$.
\qed

\lemma[lmoursuffices]%
Let $I_1, I_2, I_3$ be ideals in a Noetherian ring $A$
such that $(I_1 + I_2)^2 \subseteq I_2^2$
and let $c, c_1, \ldots, c_s$ be variables over $A$.
Then
the set of associated primes of ${A[c_1, \ldots, c_s] \over ((c_1, \ldots,
c_s) I_1 + I_2 + I_3)^n}$ equals
$$
\eqalignno{
&\left\{ (P + (c_1, \ldots, c_s)) A[c_1, \ldots, c_s]:
P \subseteq A, P + (c) \in \Ass(A[c]/(cI_1 + I_1 + I_2)^n)\right\} \cr
&\hskip 3em\cup
\left\{ P A[c_1, \ldots, c_s]:
P \in \Ass(A[c]/(cI_1 + I_1 + I_2)^n) \hbox{ and } c \not \in P\right\}. \cr
}
$$
\endb

\proof
Let $\underline c$ stand either for the ideal $(c)$
or for the ideal $(c_1, \ldots, c_s)$.
By \ref{lmpoweregen}
and using $J_1 = \underline c I_1$ and $J_2 = I_2$, $J_3 = I_3$,
we get that for all positive integers~$n$,
$(\underline c I_1 + I_1 + I_2)^n = \underline c I_1 I_3^{n-1} + (I_2 + I_3)^n$.

Define $\varphi : A[c_1, \ldots, c_s] \to A[c]$
to be the $A$-algebra homomorphism
that takes all $c_i$ to $c$ and is the identity on~$A$.
We impose the $\hbox{\tenbbold N}^s$-grading on $A[c_1, \ldots, c_s]$
with $\deg(c_i) = e_i$ (the $s$-tuple with $1$ in the $i$th position
and $0$ elsewhere)
and we define the degrees of all other variables to be $0$.
Then $\varphi$ is not a graded homomorphism,
but it is a surjective spreading as defined in \cite[Definition 2.1]{KS},
and 
$((c_1, \ldots, c_s) I_1 + I_2 + I_3)^n = (c_1, \ldots, c_s) I_1 I_3^{n-1} +
(I_2 + I_3)^n$
is a spreading of $(c I_1 + I_2 + I_3)^n = c I_1 I_3^{n-1} + (I_2 + I_3)^n$.
By \cite[Lemma 2.5]{KS},
the spreading of an irredundant primary decomposition of $c I_1 I_3^{n-1} +
(J_2+J_3)^n$
corresponds to an irredundant primary decomposition of
$(c_1, \ldots, c_s) I_1 I_3^{n-1} + (J_2+J_3)^n$;
specifically,
any associated prime of the former ideal not containing $c$
is associated to the latter ideal and does not contain any $c_i$,
and furthermore any associated prime of the former ideal that contains $c$
is spread to one unique associated prime of the latter ideal
in which the generator $c$ is replaced by the $s$ generators $c_1, \ldots, c_s$.
\qed

The last two lemmas immediately prove that the number of associated 
primes of $(\BaHH(m,r,s))^n$ is the same as
the number of associated primes of $(\BaHH(m,r,1))^n$,
via the following formalization:

\thm[thmspreadassoc]%
Set $B = \BaHH(m,r,1)$ and $B(s) = \BaHH(m,r,s)$.
For every positive integer~$n$,
the sets of associated primes of $B^n$ and of $B(s)^n$
are in one-to-one correspondence:
\ritem
Associated primes of $B^n$ not containing $c$
have the same minimal generating sets as their corresponding primes
associated to $B(s)^n$ that do not contain $c_1, \ldots, c_s$.
\ritem
Associated primes of $B^n$ that contain $c$
are of the form $P + (c)$ for some monomial prime $P$ in
variables $a_i, x_{i,j}$
and they correspond to associated primes of $B(s)^n$ of the form
$P + (c_1, \ldots, c_s)$.
\qed
\endb

\section{Lemmas}[sectlemmas]%

Throughout this section,
$B$ stands for $\BaHH(m,r,1)$ and $n$ is a positive integer.

\lemma[lm-notcontc-k1]%
Let $P$ be a prime ideal associated to $B$.
Suppose that for each $j \in [r]$
there exists $i_j \in [m]$ such that $x_{i_j, j}$ and $x_{i_j, j+1}$
are both in $P$.
Then $c \in P$.
In particular,
if $x_{i,1}$, $\ldots$, $x_{i,r} \in P$ for some $i\in [m]$,
then $c \in P$.
\endb

\proof
By definition of associated primes
there exists a monomial $w$ such that $P = (B : w)$ and hence
$w \in B : P \subseteq B : (x_{i_j,j}, x_{i_j,j+1}: j \in [r])$.
We have
$$
\eqalignno{
B : (x_{i_j,j}, x_{i_j,j+1})
&=
(B +( a_{j-1}^4x_{i_j,j-1}x_{i_j,j}, a_j^4x_{i_j,j+1}^2))
\cap (B +( a_j^4x_{i_j,j}x_{i_j,j+1}, a_{j+1}^4x_{i_j,j+2}^2))
\cr
&\subseteq
B +(a_{j-1}^4)
+(a_j^4x_{i_j,j}x_{i_j,j+1}^2,
a_j^4a_{j+1}^4x_{i_j,j+1}^2x_{i_j,j+2}^2)
\cr
&=
B +(a_{j-1}^4),
}
$$
so that
$$
w \in \bigcap_{j=1}^r (B : (x_{i_j,j}, x_{i_j,j+1}))
\subseteq \bigcap_{j=1}^r
(B +(a_{j-1}^4))
= B + (a_1^4\cdots a_r^4).
$$
In all cases,
$c$ multiplies this intersection into $B$,
proving that $c \in P$.
\qed

\lemma[lmstructurePw]%
Let $P$ be a prime ideal associated to $B^n$ containing some $x_{i,j}$.
Write $P = B^n : w$ for some monomial $w$.
\ritem
Then
$w \in a_{j-1}^4 x_{i,j-1} x_{i,j} B^{n-1} \cup a_j^4 x_{i,j+1}^2 B^{n-1}$.
\ritem
If $P$ also contains $x_{i,j+1}$,
then
$$
\eqalignno{
w &\in
\left(a_{j-1}^4 x_{i,{j-1}} x_{i,j} B^{n-1} \cup
a_j^4 x_{i,j+1}^2 B^{n-1} \right) \cap
\left(a_{j}^4 x_{i,{j}} x_{i,j+1} B^{n-1} \cup
a_{j+1}^4 x_{i,j+2}^2 B^{n-1} \right) \cr
&\hskip 5em
\cap
\left(a_{j-1}^4 x_{i,{j-1}} x_{i,j} B^{n-1} \cup
a_{j+1}^4 x_{i,j+2}^2 B^{n-1} \right).
}
$$
\endb

\proof
(1) Since $B^n : w = P$,
and by the form of the generators of $B$,
$$
w \in B^n : P \subseteq B^n : x_{i,j}
= B^n + (a_{j-1}^4 x_{i,{j-1}} x_{i,j}, a_j^4 x_{i,j+1}^2)B^{n-1}.
$$
But $w$ cannot be in $B^n$ (for otherwise $B^n : w = R$ is not a prime ideal),
and since $w$ is a monomial,
(1) follows.

(2) Now suppose that $P$ contains $x_{i,j}$ and $x_{i,j+1}$.
Then by (1),
$$
w \in
\left(a_{j-1}^4 x_{i,{j-1}} x_{i,j} B^{n-1} \cup
a_j^4 x_{i,j+1}^2 B^{n-1} \right) \cap
\left(a_{j}^4 x_{i,{j}} x_{i,j+1} B^{n-1} \cup
a_{j+1}^4 x_{i,j+2}^2 B^{n-1} \right).
$$
Suppose for contradiction that
$$
w \not \in
a_{j-1}^4 x_{i,{j-1}} x_{i,j} B^{n-1} \cup
a_{j+1}^4 x_{i,j+2}^2 B^{n-1}.
$$
Then
$w \in
a_j^4 x_{i,j+1}^2 B^{n-1}
\cap a_j^4 x_{i,j} x_{i,j+1} B^{n-1}$.
So we have proved that with $u = 0$,
$w \in
B^u (a_j^4 x_{i,j+1}^2 B^{n-1-u} \cap a_j^4 x_{i,j} x_{i,j+1} B^{n-1-u})$.
We proceed from this:
$$
\eqalignno{
w &\in
B^u (a_j^4 x_{i,j+1}^2 B^{n-1-u} \cap a_j^4 x_{i,j} x_{i,j+1} B^{n-1-u}) \cr
&=
B^u a_j^4 x_{i,j+1} \left(x_{i,j+1} B^{n-1-u} \cap x_{i,j} B^{n-1-u}\right) \cr
&=
B^u a_j^4 x_{i,j} x_{i,j+1}^2 \left(
(B^{n-1-u} : x_{i,j}) \cap (B^{n-1-u} : x_{i,j+1})
\right) \cr
&\subseteq
B^{u+1} \left(
(B^{n-1-u} +
(a_{j-1}^4 x_{i,{j-1}} x_{i,j}, a_j^4 x_{i,j+1}^2) B^{n-2-u})
\right.
\cr
&\hskip 8em
\left.
\cap (
B^{n-1-u} + (a_j^4 x_{i,j} x_{i,j+1},
a_{j+1}^4 x_{i,j+2}^2) B^{n-2-u}
)
\right). \cr
}
$$
But $w \not \in B^n \cup 
a_{j-1}^4 x_{i,{j-1}} x_{i,j} B^{n-1} \cup
a_{j+1}^4 x_{i,j+2}^2 B^{n-1}$,
so necessarily
$$
w \in
B^{u+1} (a_j^4 x_{i,j+1}^2 B^{n-2-u} \cap a_j^4 x_{i,j} x_{i,j+1} B^{n-2-u}).
$$
This is the induction step,
so when $u = n-2$,
we get that
$$
w \in
B^{n-1} ((a_j^4 x_{i,j+1}^2) \cap (a_j^4 x_{i,j} x_{i,j+1}))
= B^{n-1} (a_j^4 x_{i,j} x_{i,j+1}^2)
\subseteq B^n,
$$
which is a contradiction to $B^n :w$ being a prime ideal.
This finishes the proof of (2).
\qed

\lemma[lmstructurePgen]%
Let $P$ be a prime ideal associated to $B^n$.
Then the following properties hold:
\ritem
$a_1, \ldots, a_r \in P$.
\ritem
If $P$ does not contain $x_{i,1} x_{i,2} \cdots x_{i,r}$ for some $i$,
then $P = (a_1, \ldots, a_r)$.
\ritem[propPgentwo]%
If $P$ contains some $x_{i,j}$,
then there exists $j_0 \in \{j-1, j\}$
such that $x_{e,j_0} x_{e,j_0+1} \in P$
for all $e \in [m]$.
\ritem[propPgenfour]%
If $P$ does not contain $x_{i,j}, x_{i,j+1}$ for some $i \in [m]$
and $j \in [r]$,
then $c \not \in P$.
\endb

\proof
(1) Since $a_j^6 \in B$,
it follows that $a_j$ must be in every associated prime ideal of $B^n$.

(2) Suppose that $x_{i,1} x_{i,2} \cdots x_{i,r} \not \in P$ for some $i$.
Then $P$ is associated to
$$
B^n : (x_{i,1} \cdots x_{i,r})^\infty
= \left(B_0 + (c(a_1 a_2 \cdots a_r)^4) + (a_j^4: j \in [r])\right)^n
= (a_j^4: j \in [r])^n,
$$
whose only associated prime is $(a_1, \ldots, a_r)$.

(3) follows from \ref{lmstructurePw}~(1):
there exists $j_0 \in \{j-1, j\}$
such that the witness $w$ for $P$ is in $a_{j_0}^4 B^{n-1}$.
Hence for all $e$,
$w x_{e,j_0} x_{e,j_0+1}^2 \in B^n$.

(4) By assumption,
$P$ is associated to
$$
\eqalignno{
B^n : \left(x_{i,j} x_{i,j+1}\right)^\infty
&= \left(B + (a_j^4, a_{j-1}^4 x_{i,j-1}, a_{j+1}^4 x_{i,j+2}^2) \right)^n \cr
&= \left(
B_0 + X + (a_j^4, a_{j-1}^4 x_{i,j-1}, a_{j+1}^4 x_{i,j+2}^2)
\right)^n, \cr
}
$$
and $c$ does not appear in any generator of the last ideal.
Thus $P$ cannot contain~$c$.
\qed

\lemma[lmgij]
\ritem[lmgione]%
If $r = 2$,
then $x_{i,j} a_1^4 a_2^4 h_{i,j} \in B^2$.
\ritem[lmgitwo]%
If $r = 2$,
then $a_1^4 a_2^4 h_{i,j}^2 \in B^3$.
\ritem[lmgifive]%
$x_{i,j+2}^2 a_j^4 a_{j+1}^4 h_{i,j} \in B^2$.
\ritem[lmginine]%
$a_{j-1}^4 a_j^4 x_{i,j-1} x_{i,j} h_{i,j} \in B^2$.
\ritem[lmgithree]%
$x_{i,j} a_{j-2}^4 a_{j-1}^4 a_j^4 h_{i,j-2} h_{i,j} \in B^3$.
\ritem[lmgisix]%
$x_{i,j-1} a_{j-1}^4 a_j^4 h_{i,j}^2 \in B^3$.
\ritem[lmgifour]%
$a_{j-2}^4 a_{j-1}^4 a_j^4 h_{i,j-2} h_{i,j}^2 \in B^4$.
\endb

\proof
For \ref{lmgione} and \ref{lmgitwo} we use that $x_{i,j-1} = x_{i,j+1}$ to rewrite
$$
x_{i,j} a_1^4 a_2^4 h_{i,j} =
x_{i,j} a_1^4 a_2^4
\left(a_j^4 x_{i,j} x_{i,j+1}^2 \right)
\in (a_j^6) (a_{j-1}^4 x_{i,j-1} x_{i,j}^2) \subseteq B^2,
$$
and
$a_1^4 a_2^4 h_{i,j}^2 =
a_1^4 a_2^4
(a_j^4 x_{i,j} x_{i,j+1}^2)^2 \in
(a_j^6)^2 (a_{j-1}^4 x_{i,j-1} x_{i,j}^2) \subseteq B^3$.
Part\ref{lmgifive} follows from
$$
x_{i,j+2}^2 a_j^4 a_{j+1}^4 h_{i,j}
= x_{i,j+2}^2 a_j^4 a_{j+1}^4 (a_j^4 x_{i,j} x_{i,j+1}^2)
\in (a_j^6) (a_{j+1}^4 x_{i,j+1} x_{i,j+2}^2)
\subseteq B^2,
$$
part\ref{lmginine} from
$a_{j-1}^4 a_j^4 x_{i,j-1} x_{i,j} h_{i,j}
= a_{j-1}^4 a_j^8 x_{i,j-1} x_{i,j}^2 x_{i,j+1}^2
\in (a_j^6) (a_{j-1}^4 x_{i,j-1} x_{i,j}^2) \subseteq B^2$,
part\ref{lmgithree} from
$$
\eqalignno{
x_{i,j} a_{j-2}^4 a_{j-1}^4 a_j^4 h_{i,j-2} h_{i,j}
&=
x_{i,j} a_{j-2}^4 a_{j-1}^4 a_j^4
\left(a_{j-2}^4 x_{i,j-2} x_{i,j-1}^2 \right)
\left(a_j^4 x_{i,j} x_{i,j+1}^2 \right) \cr
&\in
(a_{j-2}^6) (a_j^6)
\left(a_{j-1}^4 x_{i,j-1} x_{i,j}^2 \right)
\subseteq B^3,
}
$$
part\ref{lmgisix} from
$x_{i,j-1} a_{j-1}^4 a_j^4 h_{i,j}^2
= x_{i,j-1} a_{j-1}^4 a_j^{12} x_{i,j}^2 x_{i,j+1}^4
\in (a_j^6)^2 (h_{i,j-1}) \subseteq B^3$,
and part\ref{lmgifour} from
$$
\eqalignno{
a_{j-2}^4 a_{j-1}^4 a_j^4 h_{i,j-2} h_{i,j}^2
&=
a_{j-2}^8 a_{j-1}^4 a_j^{12} x_{i,j-2} x_{i,j-1}^2 x_{i,j}^2 x_{i,j+1}^4 \cr
&\in
(a_{j-2}^6) (a_j^6)^2 (a_{j-1}^4 x_{i,j-1} x_{i,j}^2)
\subseteq B^4.
&\qedbox
}
$$

\cor[coraxn]%
Let $w = a_1^4 \cdots a_r^4 \left(\prod_{i,j} x_{i,j}^{v_{i,j}}\right)
\left(\prod_{i,j} h_{i,j}^{u_{i,j}}\right)$
with $v_{i,j}$, $u_{i,j}$ non-negative integers
such that $\sum_{i,j} u_{i,j} = n-1$.
\ritem
Suppose that it is possible to rewrite $w$ in the same format
but with different $v_{i,j}, u_{i,j}$.
Then $w \in B^n$.
\ritem
Suppose that $n > mr + 1$ and that $w$ multiplies
$(x_{i,j}: i \in [m], j \in [r])$ into $B^n$.
Then $w \in B^n$.
\endb

\proof
We set $A=a_1^4 \cdots a_r^4$,
$w_0 = \prod_{i,j} x_{i,j}^{v_{i,j}}$,
and $w_1= \prod_{i,j}h_{i,j}^{u_{i,j}}$.

(1)
By assumption,
there exists a positive $v_{i,j}$
such that $x_{i,j}$ gets incorporated in the rewriting of $w$
either into a new $h_{i,j-1}$ or into a new $h_{i,j}$.

Suppose that $x_{i,j}$ is incorporated into a new $h_{i,j}$.
Then $x_{i,j+1}^2$ needs to be a factor of~$w$.
This factor can come either from $u_{i,j} > 0$
or from $v_{i,j+1} + u_{i,j+1} \ge 2$.
If we use $x_{i,j+1}^2$ from $h_{i,j}^{u_{i,j}}$,
then our $x_{i,j}$ is not making a new $h_{i,j}$,
so necessarily $v_{i,j+1} + u_{i,j+1} \ge 2$.
By definition or by \ref{lmgij}\ref{lmginine} or\ref{lmgisix},
$A \cdot x_{i,j} x_{i,j+1}^2 \in (h_{i,j}) \subseteq B$,
$A \cdot x_{i,j}  x_{i,j+1} h_{i,j+1} \in B^2$,
$A \cdot x_{i,j}  h_{i,j+1}^2 \in B^3$.
This proves that $w \in B^n$.

Now suppose that $x_{i,j}$ is incorporated into a new $h_{i,j-1}$.
Then $x_{i,j-1} x_{i,j}$ needs to be a factor of~$w/x_{i,j}$.
The $x_{i,j-1}$ factor can come from $v_{i,j-1} > 0$, $u_{i,j-2} > 0$
or from $u_{i,j-1} > 0$.
We can eliminate the option $u_{i,j-1} > 0$
as it does not generate a new $h_{i,j-1}$.
Similarly, the additional factor $x_{i,j}$ can only be taken from
$v_{i,j} > 1$ or $u_{i,j} > 0$.
If $v_{i,j-1}>0$ and $v_{i,j}>1$,
then $A\cdot w_0 \in (h_{i,j-1}) \subseteq B$ and $w\in B^n$.
If $v_{i,j-1}>0$ and $u_{i,j}>0$,
then
$w \in (A\cdot x_{i,j-1}x_{i,j}h_{i,j}\cdot {w_1\over h_{i,j}}) \in B^n$
by \ref{lmgij}\ref{lmginine}.
If $u_{i,j-2}>0$ and $v_{i,j}>1$, then
$w \in (A\cdot x_{i,j}^2h_{i,j-2}\cdot {w_1\over h_{i,j-2}}) \in B^n$
by \ref{lmgij}\ref{lmgifive}.
Finally, if $u_{i,j-2}>0$ and $u_{i,j}>0$, then
$w \in (A\cdot x_{i,j}h_{i,j}h_{i,j-2}\cdot {w_1\over h_{i,j}h_{i,j-2}}) \in B^n$.

(2)
Since $n-1 > mr$,
there exists $(i,j)$ such that $u_{i,j} \ge 2$.
By assumption,
$x_{i,j}w$ and $x_{i,j+1}w$ are both in~$B^n$.
Thus for both variables,
the rewriting needs to happen as in~(1).
As in the proof of~(1),
one of the following conditions holds for $x_{i,j} w$:
\item{a)}
$v_{i,j+1} + u_{i,j+1} \ge 2$,
\item{b)}
($v_{i,j-1} > 0$ or $u_{j-2} > 0$)
and
($v_{i,j} > 0$ or $u_{i,j} > 0$);

\noindent
and one of the following conditions holds for $x_{i,j+1} w$:
\item{a')}
$v_{i,j+2} + u_{i,j+2} \ge 2$,
\item{b')}
($v_{i,j} > 0$ or $u_{i,j-1} > 0$)
and
($v_{i,j+1} > 0$ or $u_{i,j+1} > 0$).

If b) holds, then $x_{i, j-1} \mid w_0$ or $h_{i,j-2} \mid w_1$.
It follows that $w \in B^n$
since $u_{i,j} \ge 2$
and hence
$w \in (A\cdot x_{i,j-1}h_{i,j}^2\cdot {w_1 \over h_{i,j}^2}) \in B^n$
due to \ref{lmgij}\ref{lmgisix} or
$w \in (A\cdot h_{i,j-2}h_{i,j}^2\cdot {w_1 \over h_{i,j-2}h_{i,j}^2}) \in B^n$
due to \ref{lmgij}\ref{lmgifour}.
Similarly, if a') holds, then $x_{i,j+2}^2\mid w_0$ or
($x_{i,j+2}\mid w_0$ and $h_{i,j+2}\mid w_1$) or
$h_{i,j+2}^2\mid w_1$.
Since $h_{i,j}\mid w_1$,
it follows that
$w\in (A\cdot x_{i,j+2}^2h_{i,j}\cdot {w_1 \over h_{i,j}}) \in B^n$
by \ref{lmgij}\ref{lmgifive} or
$w\in (A\cdot x_{i,j+2}h_{i,j+2}h_{i,j}\cdot {w_1 \over h_{i,j}h_{i,j+2}}) \in B^n$
by \ref{lmgij}\ref{lmgithree} or
$w\in (A\cdot h_{i,j+2}^2h_{i,j}\cdot {w_1 \over h_{i,j}h_{i,j+2}^2}) \in B^n$
by \ref{lmgij}\ref{lmgifour}.

So we may assume that we have conditions a) and b').
If
$x_{i,j+1}^2\mid w_0$ and ($x_{i,j} \mid w_0$ or $h_{i,j-1} \mid w_1$),
then either
$w\in (A\cdot x_{i,j}x_{i,j+1}^2\cdot w_1\in B^n$
or
$w\in (A\cdot h_{i,j-1}x_{i,j+1}^2\cdot {w_1\over h_{i,j-1}})\in B^n$
by \ref{lmgij}\ref{lmgifive}.
If $x_{i,j+1}\mid w_0$ and $h_{i,j+1}\mid w_1$
and ($x_{i,j} \mid w_0$ or $h_{i,j-1} \mid w_1$),
then
$w\in (A\cdot x_{i,j}x_{i,j+1}h_{i,j+1}\cdot {w_1\over h_{i,j+1}})\in B^n$
by \ref{lmgij}\ref{lmginine} or
$w\in (A\cdot x_{i,j+1}h_{i,j-1}h_{i,j+1}\cdot {w_1\over h_{i,j-1} h_{i,j+1}})\in B^n$
by \ref{lmgij}\ref{lmgithree}.
Finally,
if $h_{i,j+1}^2\mid w_1$
and ($x_{i,j} \mid w_0$ or $h_{i,j-1} \mid w_1$),
then
$w\in (A\cdot x_{i,j}h_{i,j+1}^2\cdot {w_1\over h_{i,j+1}^2})\in B^n$
by \ref{lmgij}\ref{lmgisix} or
$w\in (A\cdot h_{i,j-1}h_{i,j+1}^2\cdot {w_1\over h_{i,j-1} h_{i,j+1}^2})\in B^n$
by \ref{lmgij}\ref{lmgifour}.
\qed

\lemma[lmgenwitness]%
Let $P$ be a prime ideal that contains $c$ and is associated to $B^n$.
We know that $P = B^n : w$ for some monomial $w$.
Then
\ritem[genw]%
$w = a_1^4 \cdots a_r^4 w_0$
for some $w_0 \in X^{n-1}$.

\ritem[genjjplusone]%
If $x_{i,j}$ and $x_{i,j+1}$ are in $P$
and $x_{i,j}$ is not a factor of $w$,
then $r \ge 3$ and $x_{i,j+3} \in P$.

\ritem[genjjplustwo]%
If $x_{i,j}$, $x_{i,j+1}$,  $x_{i,j+2} \in P$ and $r\ge 3$,
then $x_{i,j}$ is a factor of $w$.

\ritem[genxsquare]%
Suppose that $x_{i,j} x_{i,j+1}^2$ divides $w_0$
and that $x_{i,j+1} \in P$.
Then $x_{i,j}^2$ divides $w$.

\ritem[genxthree]%
Suppose that $x_{i,1}^2 x_{i,2}^2 \cdots x_{i,r}^2$ divides $w$.
Then $n \ge r+1$
and $w_0 \in h_{i,1} h_{i,2} \cdots h_{i,r} X^{n-1-r}$.
\endb

\proof
Since $c \in P$,
we know by \ref{lmpoweregen} that $w \in B^n : c =
a_1^4 \cdots a_r^4 X^{n-1} + (B_0 + X)^n$.
Since $w$ is a monomial not in $B^n$,
necessarily $w \in a_1^4 \cdots a_r^4 X^{n-1}$.
This proves\ref{genw}.

To simplify notation we assume in the rest of the proof that $j = 1$.

\ref{genjjplusone}
Suppose that $x_{i,1}$ does not divide $w_0$ (or $w$).
Then by \ref{lmstructurePw}~(2),
$w \in a_1^4 x_{i,2}^2 B^{n-2} \cap a_2^4 x_{i,3}^2 B^{n-2}$
and so necessarily $r \ge 3$.
This means that $w_0$ is a multiple of $x_{i,2}^2 x_{i,3}^2$.
Write $w_0 = h_{i,2}^e h_{i,3}^{e'} w'$
for some non-negative integers $e, e'$ and some $w' \in X^{n-1-e-e'}$.
We may take~$e$ to be maximal possible,
and for the maximal~$e$ we choose maximal possible~$e'$,
so that in particular $h_{i,2}$ and $h_{i,3}$ are not factors of $w'$.
By assumption also no $h_{i,1}, h_{i,r}$ appear in $w'$.
First suppose that $e = 0$.
Then by the $(x_{i,2}, x_{i,3})$-degree count,
$w_0 \in x_{i,2} x_{i,3}^2 X^{n-1}
+ x_{i,2} x_{i,3} h_{i,3} X^{n-2}
+ x_{i,2} h_{i,3}^2 X^{n-3}$,
whence by \ref{lmgij},
$w \in B^n$.
So necessarily $e \ge 1$.
Hence by \ref{lmgij}\ref{lmgifive}, $x_{i,4}^2 w \in B^n$,
so that $x_{i,4} \in P$.
This proves\ref{genjjplusone}.

We continue with the proof of\ref{genjjplustwo}. Recall that we
assume that $x_{1,i}$ does not divide $w$.
Since $x_{i,2} w \in B^n$,
necessarily in the rewriting of $x_{i,2} w$ as an element of $B^n$,
$x_{i,2}$ must combine with $x_{i,3}^2$ into a new $h_{i,2}$, i.e.,
$w/h_{i,2}^e \in a_2^4x_{i,3}^2B^{n-1-e}$.
Thus, $x_{i,2}$ is not a factor of
$w_0/h_{i,2}^e$ for otherwise $w \in B^n$.
In addition,  $x_{i,3} w \in B^n$
which implies that $x_{i,3}$ needs to recombine with $w$ into
a new element of $B$ which necessarily is $h_{i,3}$.
Thus, $w_0$ must have a factor of $x_{i,4}^2$, which comes either as
$x_{i,4}^2$, $x_{i,4} h_{i,4}$, or $h_{i,4}^2$.  But since the exponent
on $h_{i,2}$ in $w$ is at least $1$,
then by \ref{lmgij}\ref{lmgifive},\ref{lmgithree}, and\ref{lmgifour},
$w \in B^n$,
which is a contradiction, and thus proves\ref{genjjplustwo}.

\ref{genxsquare}
Suppose that $x_{i,1}^2$ does not divide $w_0$.
By assumption $x_{i,1} x_{i,2}^2$ is a factor of $w_0$,
by the $(x_{i,1}, x_{i,2})$-degree count,
$w_0 \in h_{i,1} X^{n-2} + x_{i,1} x_{i,2}^2 X^{n-1}
+ x_{i,1} x_{i,2} h_{i,2} X^{n-2}$.
If $w_0$ is in one of the last three summands,
then $w \in B^n$
by Definition or \ref{lmgij}\ref{lmginine} and\ref{lmgisix}.
So we may assume that $w_0 \in h_{i,1} X^{n-2}$.
Since $x_{i,2} w \in B^n$,
this $x_{i,2}$ must recombine with $w$ into a new $h_{i,1}$ or $h_{i,2}$,
but since there are no spare $x_{i,1}$ in $w$,
necessarily $x_{i,3}^2$ is a factor of $w_0$.
Thus, by the $x_{i,3}$-degree count,
$w_0 \in x_{i,3}^2h_{i,1}X^{n-2} +  x_{i,3}h_{i,1}h_{i,3}X^{n-3} + h_{i,1}h_{i,3}^2X^{n-4} + h_{i,1}h_{i,2}X^{n-3}$.
If $w_0$ is in the last summand, then no new $h_{i,2}$ would be form.
We can, therefore, assume that $w_0$ is in the first three summands.
But then $w \in B^n$ by \ref{lmgij}\ref{lmgifive},\ref{lmgithree}, and\ref{lmgifour}.
Thus $x_{i,1}^2$ must be a factor of~$w$.

\ref{genxthree}
Let $E$ be the largest subset of $[r]$ such that
$w_0 \in (\prod_{j \in E} h_{i,j}) X^{n-1-|E|}$.
If $E$ is empty,
then $n=1$ and
$w \in a_j^4 x_{i,j} x_{i,j+1}^2 \subseteq B = B^n$,
which is a contradiction.
Thus $E$ is not empty.
By symmetry we may assume that $h_{i,1}$ is a factor
and for contradiction we assume that $h_{i,2}$ is not a factor.
By the $(x_{i,1},x_{i,2})$-degree count,
$w_0 \in x_{i,1} h_{i,1} X^{n-2} + h_{i,1}^2 X^{n-3} + h_{i,r} h_{i,1} X^{n-3}$.
If $r = 2$,
the last summand is not possible by assumption
and the first two summands make $w$ be in $B^n$ by \ref{lmgij}\ref{lmgione} and\ref{lmgitwo},
which proves that $r \ge 3$.

If $h_{i,3}$ is also not a factor of $w_0$,
then by the $x_{i,3}$-degree count,
$w_0 \in h_{i,1} x_{i,3}^2 X^{n-2}$,
which means that $w \in B^n$ by \ref{lmgij}\ref{lmgifive}.
This proves that $h_{i,3}$ must be a factor of $w_0$,
and consequently that $E$ contains at least every other $h_{i,j}$ as $j$ varies.
Now say that $h_{i,1}, h_{i,3}$ are factors but $h_{i,2}$ is not.
By the $x_{i,3}$-degree count again,
$w_0 \in h_{i,1} x_{i,3} h_{i,3} X^{n-3} + h_{i,1} h_{i,3}^2 X^{n-3}$,
so that $w \in B^n$ by \ref{lmgij}\ref{lmgithree} and\ref{lmgifour}.
This proves\ref{genxthree}.
\qed

\section{G-good primes}[sectggood]%

The set-up is as in \ref{sectlemmas}
with $B = \BaHH(m,r,1)$ and $n$ a positive integer.
In this section
we characterize all associated primes of powers of~$B$ that are g-good.
We prove that all associated primes that contain~$c$ are g-good,
which characterizes and counts
all associated primes of powers of~$B$ that contain~$c$.
\ref{thmBHH2count} counts associated primes of any power of $\BaHH(m,2,s)$
and \ref{thmcountmax} determines the maxima
of the numbers of these associated primes.
In \ref{prop-persistence} we prove the persistence property
of associated primes of powers not containing~$c$.

We think of the $m r$ variables $x_{i,j}$
as appearing in an $m \times r$ matrix.
If a monomial prime ideal does not contain all $x_{i,j}$ in some row $i$,
then we talk about {\bf gaps},
and if a prime ideal omits some $k$ consecutive $x_{i,j}$ in a row $i$,
we refer to that as a $k$-gap of {\bf length} $k$.
Keep in mind that we identify $x_{i,r}$ with $x_{i,0}$, et cetera,
so that the gaps are counted in the round.

\defin[defgood]%
Let $P$ be a monomial prime ideal
containing $(a_1, \ldots, a_r)$.
We say that $P$ is {\bf g-good}
if it has no gaps of length 2 or larger in any of the rows.
\endb

We characterize in this section all associated prime of $B^n$
that contain~$c$ in terms of g-good primes.
The characterization enables a count, see \ref{countc}.
G-good primes also play a role for primes that do not contain~$c$;
see \refs{thmnochigherpower} and \refn{thmnocfirstpower}.

\prop[propPcfullhalf]%
Let $P$ be a g-good prime ideal containing $c$
such that for each $i \in [m]$,
the set $P \cap \{x_{i,1}, \ldots, x_{i,r}\}$
has either $r$ or exactly $r/2$ elements.
(The latter happens only if $r$ is even.)
Let $U = \{i: x_{i,1}, \ldots, x_{i,r} \in P\}$
and $V = \{(i,j): x_{i,j+1} \not \in P\}$.
\ritem
Suppose that $n = ur + v + 1$,
where $u$ and $v$ are any non-negative integers such that
$u \le |U|$ and $v \le |V|$.
Then $P$ is associated to $B^n$.
\ritem
Suppose that $n$ cannot be written as in (1).
Then $P$ is not associated to $B^n$.
\endb

\proof
Observe that $V = \{(i,j): x_{i,j} \in P\hbox{ and } i \not \in U\}$.

(1) Let $U_0$ be a subset of $U$ of cardinality $u$
and $V_0$ a subset of $V$ of cardinality $v$.
Let~$M$ be a large integer and set
$$
\eqalignno{
w_0 &= a_1^4 \cdots a_r^4
\biggl( \prod_{i \in U \setminus U_0} \prod_{j=1}^r x_{i,j} \biggr)
\biggl( \prod_{x_{i,j} \not \in P} x_{i,j}^M \biggr), \cr
w &= w_0 \biggl( \prod_{i \in U_0} \prod_{j=1}^r h_{i,j} \biggr)
\biggl( \prod_{(i,j) \in V_0} h_{i,j} \biggr). & (*)\cr
}
$$
Then $w \in B^{ur + v} = B^{n-1}$.
We will prove that $P = B^n : w$.

Since $a_j (a_1^4 \cdots a_r^4) \in (a_j^5 a_{j+1}) \in B_0 \subseteq B$
and $c (a_1^4 \cdots a_r^4) \in B_c \subseteq B$
it follows that $(a_1, \ldots, a_r, c) \subseteq B^n : w$.
Suppose that $i \in U_0$.
Then for all $j \in [r]$,
$x_{i,j} \in B^n : w$
by \ref{lmgij}\ref{lmgione} in case $r = 2$
and by \ref{lmgij}\ref{lmgithree} in case $r > 2$.
If $i \in U \setminus U_0$,
then
$x_{i,j} w \in x_{i,j} (a_{j-1}^4 x_{i,j-1} x_{i,j}) X^{n-1} \subseteq
X^n \subseteq B^n$.
Thus $x_{i,j} \in B^n : w$ for all $i \in U$ and all $j \in [r]$.
If $(i,j) \in V$,
then $x_{i,j+1} \not \in P$,
so that $x_{i,j} w_0 \in (a_j^4 x_{i,j} x_{i,j+1}^2) \subseteq X$
and thus $x_{i,j} \in B^n : w$.
This proves that $P \subseteq B^n : w$.

To prove that $P = B^n : w$
it remains to show that every $x_{i,j} \not \in P$
is a non-zerodivisor modulo $B^n$.
By possibly taking $M$ even larger
it suffices to prove that $w \not \in B^n$.
In the given form $w$ is an element of~$B^{n-1}$.
Any rewriting of $w$ to make it an element of $B^n$
has to involve the variables $x_{i,j}$
whose exponents are at least two.
The only such $x_{i,j}$ are those not in $P$ and those with $i \in U_0$.
If $x_{i,j} \not \in P$,
then $x_{i,j+1}^2$ is not a factor of $w$
and either $x_{i,j-1}$ is not a factor of $w$
or else
$x_{i,j-1}$ is a factor of $w$ but tied up in $h_{i,j-1}$.
Thus there is no possible way of using the rewriting with $x_{i,j}$ not in~$P$.
If $i \in U_0$,
then this $i$ contributes to $w$ the factor
$a_1^4 \cdots a_r^4 x_{i,1}^3 \cdots x_{i,r}^3 \in B^r$,
and
there is no possible way of rewriting this part to put $w$ into $B^n$.
Thus $P = B^n : w$.
It follows that $P$ is associated to $B^n$.
This finishes the proof of (1).

(2)
Let $P$ be associated to $B^n$
and suppose for contradiction that $n$ cannot be written as in (1).
Then in particular $n > 1$.
By \ref{lmgenwitness}\ref{genw},
$P = B^n : w$,
where $w$ is a monomial of the form $a_1^4 \cdots a_r^4 w_0$
for some $w_0 \in X^{n-1}$.
The product of all $x_{i,j}$ with $i \in U$ divides $w_0$
by \ref{lmgenwitness}\ref{genjjplusone} in case $r = 2$
and by \ref{lmgenwitness}\ref{genjjplustwo} in case $r > 2$.
Let $U_0$ be the set of all $i \in U$ such that
$h_{i,j}$ is a factor of $w_0$ for some $j \in [r]$.
Let $i \in U_0$.
Then by \ref{lmgenwitness}\ref{genxsquare},
$x_{i,j}^2$ divides $w_0$.
Since $x_{i,j-1}$ also divides $w_0$,
then again by \ref{lmgenwitness}\ref{genxsquare},
$x_{i,j-1}^2$ divides $w_0$.
By continuing in this way we get that
$\prod_{j=1}^r x_{i,j}^2$ divides $w_0$.
Hence by \ref{lmgenwitness}\ref{genxthree},
$\prod_{j=1}^r h_{i,j}$ divides $w_0$.
By \ref{lmgij}\ref{lmgitwo} and\ref{lmgifour},
$h_{i,j}^2$ is not a factor of $w$ for all such $i, j$.
Thus $n - 1 \ge |U_0| r$.
Let $v = n - 1 - |U_0| r$.
We just proved that $w_0$ is a product of the $|U_0| r$
factors $h_{i,j}$ with $(i,j) \in U_0 \times [r]$
and $v$ factors $h_{i,j}$ with $i \not \in U$.
By \ref{lmgij}\ref{lmgifive},
necessarily for any such latter factor we have
$x_{i,j+2}^2 \in B^n : w = P$.
Since $P$ is g-good and the $i$th row has $r/2$ elements,
necessarily $x_{i,j} \in P$ and $x_{i,j-1}$ is not in $P$.
Then by \ref{lmgij}\ref{lmgisix},
the squares of these $h_{i,j}$ do not divide $w$.
Thus these $v$ factors are all distinct,
which means that $n$ must be written as in (1).
\qed

\thm[thmcfullhalf]%
We consider the set $S$ of all g-good prime ideals~$P$ containing~$c$
for which in each row of the matrix $[x_{i,j}]$,
$P$ contains either $r$ or $r/2$ elements.
\ritem
If $r$ is odd,
then the maximal ideal is the only such prime ideal,
and it is associated to $B^n$
if and only if $n = ur+1$ for some $u \le m$.
\ritem
If $r$ is even,
then $S$ contains $3^m$ prime ideals.
\riitem{(a)}
For each $i \in \{0,\ldots, m\}$
there exist $2^i {m \choose i}$ prime ideals in $S$
of height $(m+1) r - i {r \over 2} + 1$,
and these are associated to $B^n$ exactly when $n$
equals $u r + v + 1$ with $u \in \{0, \ldots, m-i\}$
and $v \in \{0, \ldots, i{r \over 2}\}$.
\riitem{(b)}
The number $h(m,r,n)$ of elements of $S$
that are associated to $B^n$ equals
$$
\sum_{i = 0}^m 2^i {m \choose i} \delta_{(n-1)/r - i/2 \le \min\{q, m-i\}},
$$
where
$q = \lfloor \frac{n - 1}{r}\rfloor$.
For all $n > 1 + rm$, $h(m,r,n) = 0$.
\endb

\proof
(1) is an immediate corollary of \ref{propPcfullhalf}.

To prove (2),
observe that for $i \in [m]$,
one of three things happen for $P \in S$:
$P$ contains the full $i$th row of $[x_{i,j}]$,
$P$ contains $x_{i,j}$ with $j$ odd,
and
$P$ contains $x_{i,j}$ with $j$ even.
Thus the count of elements of $S$ is $3^m$.
For each $i \in \{0, 1, \ldots, m\}$,
there are ${m \choose i}$ possibilities
where exactly $m-i$ of the rows are fully in~$P$,
and the remaining $i$ rows have two options.
All these prime ideals contain also $a_1, \ldots, a_r, c$,
so that their height is
$r + 1 + (m-i)r + i{r \over 2}$
$= (m+1) r - i {r \over 2} + 1$.

According to \ref{propPcfullhalf},
$P$ is associated to $B^n$ if and only if
there exist integers $u \in \{0, \ldots, m-i\}$
and $v \in \{0, \ldots, i\frac{r}{2}\}$
such that $n-1 = ur+v$.
The rest of (2)(a) is an immediate corollary of \ref{propPcfullhalf}.

For (2)(b) we need to account which $n$ are possible.
Note that $n-1 = ur+v \le (m-i)r + ir/2 \le mr$.
Thus $h(m,r,n) = 0$ if $n - 1 > mr$.
The possible $u$ are $0, 1, \ldots, m-i$,
if simultaneously $0 \le v = n-1-ur \le ir/2$.
Another way of recording this is with
$\max\{0,(n-1)/r - i/2\} \le u \le \min\{q, m-i\}$.
The assertion in (2)(b) follows because $\min\{q, m-i\}\ge 0$.
\qed

\thm[thmotherc]
Let $P$ be a g-good monomial prime ideal containing $c$.
Suppose that there exists $(i_0, j_0) \in [m] \times [r]$
such that $x_{i_0, j_0-1}, x_{i_0, j_0} \in P$ and
$x_{i_0, j_0+1} \not \in P$.
Then $P$ is associated to $B^n$ for all $n \ge 1$.
\endb

\proof
The assumption on $i_0, j_0$ forces $r \ge 3$.
Let $T_1$ be the set of all $x_{i,j}$ not in $P$
and $T_2$ the set of all $x_{i,j} \in P$ such that $x_{i,j+1} \in P$.
In case $n \ge 3$
we correct~$T_2$ to not include $x_{i_0, j_0-1}$.
For any large integer $M$ we set
$$
w_0 = a_1^4 \cdots a_r^4
\left( \prod_{t \in T_1} t^M \right)
\left( \prod_{t \in T_2} t \right),
\hskip 1em
w = w_0 h_{i_0, j_0}^{n-1}.
$$
Since $(a_1, \ldots, a_r, c) a_1^4 \cdots a_r^4 \in B$,
it follows that $(a_1, \ldots, a_r, c) \subseteq B^n : w$.
If $x_{i,j} \in P$ and $x_{i,j+1} \not \in P$,
then $x_{i,j} w \in B^n$ since $a_j^4 x_{i,j+1}^2$ is a factor of $w_0$
and $h_{i,j} = a_j^4 x_{i,j} x_{i,j+1}^2 \in B$.
If $n \ge 3$,
then
$x_{i_0, j_0-1} w
\in (h_{i_0, j_0 - 1}) (a_{j_0}^{12}) (h_{i_0, j_0}^{n-3})
\subseteq B^{1 + 2 + n - 3} = B^n$.
In all other cases,
if $x_{i,j}, x_{i,j+1} \in P$,
then $x_{i,j} w \in B^n$
since $a_{j-1}^4 x_{i,j-1} x_{i,j}$ is a factor of $w_0$
and $h_{i,j-1} = a_{j-1}^4 x_{i,j-1} x_{i,j}^2 \in B$.
This proves that $P \subseteq B^n : w$.

We next prove that $B^n : w \subseteq P$,
i.e., that no power of a variable in $T_1$ is in $B^n : w$.
By possibly taking~$M$ larger it suffices to prove that $w \not \in B^n$.
In the given form $w$ is an element of~$B^{n-1}$.
Any rewriting of $w$ to make it an element of $B^n$
has to involve the variables $x_{i,j}$
whose exponents are at least two.
The only such variables are those in $T_1$
and additionally $x_{i_0, j_0}$ if $n-1 \ge 2$.
By the g-goodness assumption,
the variables in $T_1$ do not have consecutive second indices
and $T_2$ does not contain suitable ``predecessors''
to form a new $h_{i,j-1}$ with a variable $x_{i,j}\in T_1$.
So necessarily $n \ge 3$,
but then $x_{i_0, j_0-1}$ is not a factor of $w$
so it is not possible to recombine $x_{i_0, j_0}^2$ with that missing factor
and no other rewriting is possible.
Thus $w \not \in B^n$.

Thus $P = B^n : w$ so that $P$ is associated to $B^n$.
\qed

\lemma[lmLucas]%
The number of g-good primes
(either all containing $c$ or none containing~$c$)
equals the Lucas number $L_r^m$
(with $L_1 = 1$, $L_2 = 3$, $L_{r+2} = L_{r+1} + L_r$).
\endb

\proof
Note that the number of g-good primes of either type is equal to $L_r^m$
where $L_r$ is their number for the case $m=1$.
We will ignore containments of $a_1, \ldots, a_r$ in this proof.

In case $r = 1$,
the only g-good prime contains $x_{1,1}$,
so $L_1$ is~$1$.
The g-good options in case $r = 2$
are $(x_{1,1})$, $(x_{1,2})$, and $(x_{1,1}, x_{1,2})$,
so $L_2 = 3$.
Let $U_r$ be the number of g-good primes that contain $x_{1,1}$ and $x_{1,r}$,
and for $r > 1$ let $V_r$ be the number of g-good primes that
contain $x_{1,1}$ and not $x_{1,r}$, and let $\overline V_r$ be the number of
g-good primes that contain $x_{1,r}$ and not $x_{1,1}$.
Clearly $\overline V_r
= V_r$, $V_r = U_{r-1}$, and $L_r = U_r + V_r + \overline V_r = U_r +
2 V_r$.
But $U_{r+1} = U_r + U_{r-1}$
(depending on whether $r-1$ is or is not in the subset),
$U_1 = 1$,
$U_2 = 1$,
and $U_3 = 2$,
which says that $U_1, U_2, U_3, \ldots$ are the usual Fibonacci numbers,
and so $V_2, V_3, V_4, \ldots$ are also the usual Fibonacci numbers.
Then
$$
L_{r+1} + L_r
= U_{r+1} + 2 U_r + U_r + 2 U_{r-1}
= (U_{r+1} + U_r) + 2 (U_r + U_{r-1})
= U_{r+2} + 2 U_{r+1} = L_{r+2},
$$
and so these numbers are the Lucas numbers.
\qed

\theorem[countc]
The number of prime ideals associated to $\BaHH(m,r,1)^n$
that contain $c$ is equal to
$$
\cases{
L_r^m - 3^m + h(m,r,n), & if $r$ even; \cr
L_r^m, & if $r$ odd, $n \equiv 1 \mod r$ and $n \le rm+1$; \cr
L_r^m -1, & otherwise, \cr
}
$$
where $L_r$ is the $r$th Lucas number with $L_1=1$ and $L_2=3$
and $h(m,n,r)$ refers to the number in \ref{thmcfullhalf}~(b).
\endb

\proof
By \ref{lmstructurePgen}\refn{propPgenfour},
every prime ideal associated to $B^n$
that contains $c$ must be g-good.

Assume first that $r$ is odd.
In this case,
according to \ref{thmcfullhalf}~(1),
the maximal ideal
is associated if and only if $n = ur + 1$ for some integer $u \le m$.
Any other one of the $L_r^m$ possible prime ideals
satisfies condition of \ref{thmotherc}
and is thus associated to $B^n$ for all $n$.
This proves the theorem for odd $r$ by \ref{lmLucas}.

Now, assume that $r$ is even.
Of the $L_r^m$ possible prime ideals
as accounted for by \ref{lmLucas},
those for which some row in the matrix $[x_{i,j}]$ is neither half-full nor full
are covered by \ref{thmotherc}
and are thus associated to all powers of $B$.
It remains to count those prime ideals associated to $B^n$
for which each row in $[x_{i,j}]$ is either half full or full.
According to \ref{thmcfullhalf}~(2),
there are $3^m$ prime ideals with only full and half-full levels,
of which $h(m,r,n)$ are associated to $B^n$.
The theorem follows.
\qed

\example
The following tables
of numbers of associated primes of $\BaHH(m,r,1)^n$
that contain~$c$
are taken from \ref{countc}
and agree with the calculations%
\footnote *{The program code associated with this paper is
available as ancillary file on the arXiv page of this paper.}
by Macaulay2~\cite{GS} and Magma~\cite{Magma}
of associated primes for low values of $n$.

$r =2$
\vskip 1ex

\def\vr{\smash{\vrule height 2.1ex depth 1.3ex width .05em}}
\halign{\vadjust{\vskip2pt\hrule}\vr\ \hfil$#$\hfil\ \vr
&& \hskip3pt\hfil$#$ \hfil \vr \cr
\noalign{\hrule\vskip0.2ex}
m\backslash n & 1 & 2 & 3 & 4 & 5 & 6 & 7 & 8 & 9 & 10 & 11 & 12 & 13 & 14 \cr
\noalign{\vskip 1pt\hrule\vskip 2pt}
1 & 3 & 2 & 1 & 0 & 0 & 0 & & & & & & & & \cr
2 & 9 & 8 & 9 & 4 & 1 & 0 & & & & & & & & \cr
\noalign{\gdef\vr{\smash{\vrule height 2.8ex depth 0.5ex width .05em}}}
3 & 27 & 26 & 27 & 26 & 19 & 6 & 1 & 0 & & & & & & \cr
4 & 81 & 80 & 81 & 80 & 81 & 64 & 33 & 8 & 1 & 0& & & & \cr
5 & 243 &242 & 243 & 242 & 243 & 242 & 211 & 130 & 51& 10& 1 & 0 & & \cr
6 & 729& 728& 729& 728& 729& 728& 729& 664& 473& 232& 73& 12& 1 & 0 \cr
}

$r =4$
\vskip 1ex

\def\vr{\smash{\vrule height 2.1ex depth 1.3ex width .05em}}
\halign{\vadjust{\vskip2pt\hrule}\vr\hfil$#$\hfil\vr
&& \hskip1pt\hfil$#$ \hfil \vr \cr
\noalign{\hrule\vskip0.2ex}
m\backslash n & 1 & 2 & 3 & 4 & 5 & 6 & 7 & 8 & 9 & 10 & 11 & 12 \cr
\noalign{\vskip 1pt\hrule\vskip 2pt}
1 &7 & 6 & 6 & 4 & 5 & 4 & 4 & 4 & 4 & 4 & 4 & 4 \cr
2 & 49 & 48 & 48 & 44 & 49 & 44 & 44 & 40 & 41 & 40 & 40 & 40 \cr
3 & 343 & 342 & 342 & 336 & 343 & 342 & 342 & 328 & 335 & 322 & 322 & 316 \cr
4 &2401 & 2400 & 2400 & 2392 & 2401 & 2400 & 2400 & 2392 & 2401 & 2384 & 2384 & 2344 \cr
\noalign{\gdef\vr{\smash{\vrule height 2.8ex depth 0.5ex width .05em}}}
5 &16807 & 16806 & 16806 & 16796 & 16807 & 16806 & 16806 & 16796 & 16807 & 16806 & 16806 & 16764 \cr
}

\vskip 2ex
We have finished a characterization of all associated primes of $B^n$
that contain~$c$.

In contrast,
we do not have a complete characterization of the prime ideals
associated to $B^n$ that do not contain~$c$.
Of these,
we understand the g-good ones well:
by \ref{thmnochigherpower},
the number of such is $L_r^n$ if $n \ge 2$,
but the count is smaller for $n = 1$ by \ref{thmnocfirstpower}.

\thm[thmnochigherpower]%
Let $P$ be a g-good monomial prime ideal that does not contain~$c$
and let $n \ge 2$.
Then $P$ is associated to~$B^n$.
The number of such primes is~$L_r^n$.
\endb

\proof
Set
$$
e_n = \cases{
5n-5, & if  $n = 2, 3, 4$; \cr
6n-9, & if $n \ge 4$. \cr
}
$$
Note that
$$
\eqalignno{
a_1^{e_n+3} &=
\cases{
a_1^{5n-2} \in (a_1^6)^{n-1} \subseteq B^{n-1}, & if $n = 2$, $3$, $4$; \cr
a_1^{6n-6} \in (a_1^6)^{n-1} \subseteq B^{n-1}, & if $n \ge 4$, \cr
},
\cr
a_1^{e_n} a_2^4 &=
\cases{
a_1^{5(n-1)} a_2^4 \in (a_1^5a_2)^{n-1} \subseteq B^{n-1}, & if $n = 2$, $3$, $4$; \cr
a_1^{6n-9} a_2^4 \in (a_1^5 a_2)^4 (a_1^6)^{n-5} \subseteq B^{n-1}, & if $n \ge 5$, \cr
}\cr
a_1^{e_n + 5} a_2^4
&= \cases{
a_1^{5n} a_2^4 \in (a_1^5a_2)^n \subseteq B^n, & if $n = 2$, $3$, $4$; \cr
a_1^{6n-4} a_2^4 \in (a_1^5 a_2)^4 (a_1^6)^{n-4} \subseteq B^n, & if $n \ge 4$. \cr
}
\cr
}
$$
With $w_0 = a_1^{e_n} a_1^4 \cdots a_r^4$,
we have that
$$
\eqalignno{
a_1 w_0 &\in (a_1^{e_n+5} a_2^4) \subseteq B^n, \cr
a_j w_0 &\in (a_j^5 a_{j+1}) a_1^{e_n+3} \subseteq B^n
\hbox{ if $j \in \{2, \ldots, r\}$}, \cr
}
$$
so that $(a_1, \ldots, a_r) w_0 \subseteq B^n$.

Let $X_P$ be the product of powers of the $x_{i,j}$,
where
$$
\hbox{the exponent of $x_{i,j}$ in $X_P$}
=
\cases{
2, & \hbox{ if $x_{i,j} \not \in P$}; \cr
1, & \hbox{ if $x_{i,j} \in P$ and $x_{i,j+1} \in P$}; \cr
0, & \hbox{ if $x_{i,j} \in P$ and $x_{i,j+1} \not \in P$}. \cr
}
$$
Set $w = c w_0 X_P$.
We will prove that $P = B^n : w$.
We have established that $(a_1, \ldots, a_r) w \subseteq B^n$.
Now let $x_{i,j} \in P$.
If $x_{i,j+1} \not \in P$,
then $x_{i,j+1}^2$ divides $X_P$,
and so
$$
x_{i,j} w \in (a_j^4 x_{i,j} x_{i,j+1}^2) \cdot
\cases{
a_1^{e_n} a_2^4, & if $j = 1$; \cr
a_1^{e_n+3}, & if $j \not = 1$, \cr
}
$$
which is in $B^n$.
If instead $x_{i,j+1} \in P$,
then $x_{i,j-1}x_{i,j}$ divides $X_P$,
and so
$$
x_{i,j} w \in (a_{j-1}^4 x_{i,j-1} x_{i,j}^2) \cdot
\cases{
a_1^{e_n} a_2^4, & if $j = 2$; \cr
a_1^{e_n+3}, & if $j \not = 2$, \cr
}
$$
which is again in $B^n$.
This proves that $P w \subseteq B^n$.

To prove that $B^n : w \subseteq P$,
we need to prove that $z w \not \in B^n$,
where $z$ is a high power of a product of $c$
and all the $x_{i,j}$ that are not in~$P$.
Any rewriting of $zw$ as an $n$-fold product of elements in $B$
cannot use any generator of $\underline a$-degree~4
because $X_P$ contains no factors of the form $x_{i,j} x_{i,j+1}^2$.
Thus the $n$ factors in this rewriting are taken from the following list:
$a_1^6, a_1^5 a_2, c a_1^4 \cdots a_r^4$.
If the latter factor appears,
then $zw$ would have to be a multiple of $(a_1^6)^{n-1} c a_1^4 \cdots a_r^4$,
but the $a_1$-degree is then too high.
Thus the only possible factors are
$a_1^6$ and $a_1^5 a_2$.
It is easy to see that this is not possible if $n = 2, 3$,
and for $n \ge 4$,
the total degree $e_n + 8 = 6n-1$ of $a_1$ and $a_2$ in $zw$
would have to be at least $6n$, which is a contradiction.
This finishes the proof that $B^n : w = P$,
so that $P$ is associated to $B^n$.

The number of such primes was determined in \ref{lmLucas}.
\qed

\thm[thmnocfirstpower]%
Let $P$ be a g-good monomial prime ideal that does not contain~$c$.
Then $P$ is associated to $B$ if and only if
there exists $j_0 \in [r]$ such that
for all $i \in [m]$,
either $x_{i,j_0} \not \in P$ or $x_{i,j_0+1} \not \in P$.

When $r = 2$, the number of such $P$ is exactly $2^m$.
\endb

\proof
By \ref{lm-notcontc-k1},
if $P$ is associated then such a $j_0$ must exist.
Now suppose that $j_0$ exists.
By possibly replacing $j_0$ with $j_0 + 1$ we may assume that
there exists $i \in [m]$ such that $x_{i,j_0} \not \in P$.
By re-indexing we may assume that $j_0 = 1$.

Let $X_P$ be the product of various powers of the $x_{i,j}$,
where
$$
\hbox{the exponent of $x_{i,j}$ in $X_P$}
=
\cases{
2, & \hbox{ if $x_{i,j} \not \in P$}; \cr
1, & \hbox{ if $x_{i,j} \in P$ and $j = r$}; \cr
1, & \hbox{ if $x_{i,j} \in P$, $j \not = r$ and $x_{i,j+1} \in P$}; \cr
0, & \hbox{ if $x_{i,j} \in P$, $j \not = r$ and $x_{i,j+1} \not \in P$}. \cr
}
$$
and set
$w = c a_r^3 \cdot \left( \prod_{j=1}^{r-1} a_j^4 \right) X_P$.
We will prove that $P = B : w$.

We have that $a_r w = (ca_1^4 a_2^4\cdots a_r^4) \subseteq B$,
and for $j \in \{1, \ldots, r-1\}$,
we have that $a_j w \in (a_j^5 a_{j+1}) \subseteq B$.
This proves that $(a_1, \ldots, a_r) w \subseteq B$.

Now let $x_{i,j} \in P$.
We need to prove that $x_{i,j} w \in B$.
If $j = r$,
then $x_{i,r-1} x_{i,r}$ divides~$X_P$ and $a_{r-1}^4$ divides~$w$.
Hence,
$x_{i,r} w \in (a_{r-1}^4 x_{i,r-1} x_{i,r}^2) \subseteq B$.
If $j = 1$,
by the definition of $j_0$,
$x_{i,2}$ is not in~$P$.
Thus $a_1^4 x_{i,2}^2$ is a factor of $w$,
so that
$x_{i,1} w \in (a_1^4 x_{i,1} x_{i,2}^2) \subseteq B$.
Now let $j \in \{2, \ldots, r-1\}$.
We need to prove that $x_{i,j} w \in B$.
So necessarily $r > 2$.
Then $X_P$ is a multiple of $x_{i,j+1}^2$ if $x_{i,j+1} \not \in P$,
or else it is a multiple of $x_{i,j-1} x_{i,j}$,
so that
$$
x_{i,j} w \in \cases{
(a_j^4 x_{i,j} x_{i,j+1}^2) \subseteq B,
& if $x_{i,j+1} \not \in P$; \cr
(a_{j-1}^4 x_{i,j-1} x_{i,j}^2) \subseteq B,
& if $x_{i,j+1} \in P$. \cr
}
$$
This finishes the proof that $P \subseteq B : w$.

It remains to prove that $B : w \subseteq P$.
It suffices to prove that $z w \not \in B$,
where $z$ is a high power of a product of $c$
with all $x_{i,j}$ that are not in~$P$.
Since the $a_r$-degree of $w$ is $3$,
the rewriting of $zw$ as an element of $B$
would not use the one generator of $B$ that involves $c$.
So $c$ plays no role in this rewriting.
If $x_{i,j} \not \in P$,
then
both $x_{i,j-1}$ and $x_{i,j+1}$ must be in $P$.
Thus $x_{i,j+1}^2$ is not a factor of $zw$
and $x_{i,j-1}$ is a factor of $zw$ exactly if $j-1 = r$.
In that case,
$a_{j-1}^4$ is not a factor of $zw$,
which means that no rewriting of $zw$ as an element of~$B$
can use $a_{j-1}^4 x_{i,j-1} x_{i,j}^2$
or $a_j^4 x_{i,j} x_{i,j+1}^2$.
Also,
no factor of this form already appears in~$w$.
But then by the consideration of exponents of the $a_j$ in $w$,
$zw \not \in B$.
\qed

We have handled all the g-good prime ideals associated to powers of~$B$.
There are further associated primes that do not contain~$c$.
In the rest of this section we prove their persistence property
and we completely describe and enumerate them in case $r = 2$.
Persistence definitely fails on associated primes that do contain~$c$.

\proposition[prop-persistence]
\thmparens{Persistence of associated primes}
Let $P$ be a prime ideal associated to $B^n$ that does not contain~$c$.
Then $P$ is associated to $B^{n+1}$.
\endb

\proof
If $P$ is g-good,
then
$P$ is associated to all $B^{n+1}$
by \ref{thmnochigherpower}.

So we may assume that there exists $i \in [m]$ and $j \in [r]$
such that $x_{i,j}, x_{i,j+1} \not \in P$.
Then $P$ is associated to $B^n$ if and only if
it is associated to $B^n$ after inverting $x_{i,j} x_{i,j+1}$.
But then $a_j^4$ is a minimal generator of $B$,
and the only other minimal generator of $B$ (after this inversion)
in which $a_j$ appears is
$a_{j-1}^5 a_j$.
Write $P = B^n : w$ for some monomial $w$.
Then
$$
\eqalignno{
B^{n+1} : a_j^4 w
&= (B^{n+1} : a_j^4) : w \cr
&= (B^n + a_{j-1}^5 B^n
+ a_{j-1}^{10} B^{n-1}
+ a_{j-1}^{15} B^{n-2}
+ a_{j-1}^{20} B^{n-3}) : w \cr
&= B^n : w = P, \cr
}
$$
so that $P$ is associated to $B^{n+1}$ as well.
\qed

\thm[thmBHH2count]%
Let $m$, $s\ge 1$. The set of associated primes of $\BaHH(m,2,s)^n$
is the union $\{(a_1, a_2)\} \cup \bigcup_{i=1}^m Q_c^{(i)} \cup Q_1 \cup Q_2$ where
$$
\eqalignno{
Q_c^{(i)} &= \{P \mid n=2u+v+1 \text{ with } 0\le u \le m-i, 0\le v\le i \text{ and} \cr
          &\quad\quad\quad\; \text{$P$ has $i$ half-full and $m-i$ full rows}\},  \cr
Q_1&=\{P \text{ g-good} \mid c\notin P, n=1, \exists j_0\in[r]\, \forall i\in[m]\colon x_{i,j_0}\notin P \text{ or } x_{i,j_0+1}\notin P\} \cr
Q_2&=\{P \text{ g-good} \mid c\notin P, n\ge 2\}
}
$$

\noindent
The number of associated primes of $\BaHH(m,2,s)^n$ is equal to
$$
(3-\delta_{1=n})^m +
\sum_{i = 0}^m 2^i {m \choose i} \delta_{(n-1-i)/2 \le \min\{q, m-i\}} + 1.
$$
\endb

\proof
By \ref{thmspreadassoc},
the set of associated primes of $\BaHH(m,2,s)^n$ is equal to
the set of associated primes of $\BaHH(m,2,1)^n$.
By \refs{countc} and \refn{thmcfullhalf},
the set of associated primes of $\BaHH(m,2,1)^n$ that contain $c$ equals
$\bigcup_{i=1}^m Q_c^{(i)}$ and its cardinality is
$h(m,2,n) = 
\sum_{i = 0}^m 2^i {m \choose i} \delta_{(n-1-i)/2 \le \min\{q, m-i\}}$,
where
$q = \lfloor \frac{n - 1}{2}\rfloor$.

Let $P$ be associated to $\BaHH(m,2,1)^n$ and not contain~$c$.
If $P$ does not contain $x_{i,1}x_{i,2}$,
then $P$ is associated to $B^n$
if and only if it is associated to
$B^n : (x_{i,1} x_{i,2})^\infty = (a_1^4, a_2^4)^n$,
in which case $P$ must be equal to $(a_1, a_2)$,
which is minimal over $B$ and hence associated to all the powers of~$B$.

Thus it remains to consider the associated primes $P$ not containing $c$
that contain $x_{i,1} x_{i,2}$ for all $i \in [m]$.
Then $P$ must be g-good,
and in $Q_1$ if $n = 1$ by \ref{thmnocfirstpower} (which has cardinality $2^m$)
and in $Q_2$ if $n \ge 2$ by \ref{thmnochigherpower} (which has cardinality $3^m$).

The assertion follows.
\qed

\remark[rmknewcomb]%
The number of associated primes of $\BaHH(m,2,s)^n$ can also be written as
$$
(3-\delta_{1=n})^m +
\left( \sum_{\ell=0}^m \sum_{t=b(\ell)}^m { {m}\choose{\ell}} {{\ell}\choose{\ell+t-m} } \right) +
\cases{0, & if $n \le 2m$ and $n$ is even; \cr
1, & otherwise. \cr }
$$
where $b(\ell) = \max\{n-1-\ell,m-\ell\}$.
Namely,
the first summand in the display plus 1 is the number of associated primes not
containing~$c$.
The maximal ideal is associated to $\BaHH(m,2,1)^n$
if and only if $n \le 2m+1$ and $n$ is odd.
These two counts account for the first and the last summand in the display.
It remains to count the non-maximal associated primes~$P$ that contain~$c$.
We know that for all $i \in [m]$,
$x_{i,1}x_{i,2} \in P$.
Let $\ell$ be the number of $x_{i,1}$ in $P$
and let $t$ be the number of $x_{i,2}$ in $P$.
Necessarily $\ell+t \ge m$.
Also,
$\ell + t$ should be at least $2u + v = n - 1$
as in the notation of \ref{thmcfullhalf}.
There are ${m \choose \ell}$ ways of choosing $\ell$ of the variables $x_{i,1}$,
after which for the remaining $m - \ell$ rows in the matrix $[x_{i,j}]$,
the elements $x_{i,2}$ must be in~$P$.
This leaves $t - (m-\ell)$ variables $x_{i,2}$ to be chosen
from the $\ell$ rows with the $x_{i,1}$.
This justifies the middle summand in the display.

We just proved the following combinatorial identity:

$$
\eqalignno{
&\sum_{i = 0}^m 2^i {m \choose i} \delta_{n-1 \le \min\{2q+i, 2m-i\}}
=\cr
&\sum_{\ell=0}^m \sum_{t=b(\ell)}^m { {m}\choose{\ell}} {{\ell}\choose{\ell+t-m} }
+
\cases{-1, & if $n \le 2m$ and $n$ is even; \cr
0, & otherwise,\cr }
}
$$
where $b(\ell) = \max\{n-1-\ell,m-\ell\}$.

\thm[thmcountmax]
For $m \ge 1$,
the function $\phi$ taking
$n \mapsto \#\Ass(R/\BaHH(m,2,s)^n)$
has exactly $\left\lceil {m-1 \over 2} \right\rceil$ local maxima.
The local maxima occur
at $n = 3, 5, \ldots, 2 \left\lceil {m-1 \over 2} \right\rceil + 1$,
and they are all equal to the global maximum
$2 \cdot 3^m + 1$.
\endb

\proof
We refer to the three summands in the display in \ref{rmknewcomb}
as $\phi_1, \phi_2, \phi_3$ (in the given order).
Observe that $\phi_1$ is constant for $n \ge 2$,
that $\phi_2$ is zero for all $n \ge 2m+2$,
and that $\phi_3$ is constant for $n \ge 2m+1$.
Thus $\phi$ is constant for $n \ge 2m+2$.

In the range $n = 1, \ldots, m+1$,
$\phi_2$ equals
$$
\sum_{\ell=0}^m \sum_{t=m-\ell}^m { {m}\choose{\ell}} {{\ell}\choose{\ell+t-m} }
= (1 + 1 + 1)^m = 3^m,
$$
after which it strictly decreases to $0$ at $n = 2m+2$.
Thus $\phi(1) = 2^m + 3^m + 1 < 2 \cdot 3^m = \phi(2)
< 2 \cdot 3^m + 1 = \phi(3)$,
and this is equal to $\phi(n)$ for all odd $n \in \{3, \ldots, m+1\}$.
In other words,
$\phi(n) = \phi(3)$
for all $n = 3, 5, \ldots, 2 \left\lceil {m-1 \over 2} \right\rceil + 1$.
This value is strictly larger than
$2 \cdot 3^m = \phi(4) = \phi(6) = \cdots
= \phi(2 \left\lceil {m-1 \over 2} \right\rceil)$,
and is also strictly larger than
$\phi(2 \left\lceil {m-1 \over 2} \right\rceil + 2)$.

Furthermore,
for $n \in \{m+1, m+2, \ldots, 2m\}$,
$$
\phi_2(n) - \phi_2(n+1)
= \sum_{\ell=0}^m { {m}\choose{\ell}} {{\ell}\choose{\ell+(n-1- \ell)-m} }
= \sum_{\ell=n-1-m}^m { {m}\choose{\ell}} {{\ell}\choose{n-1-m}}
\ge 2.
$$
Thus $\phi(n) > \phi(n+1)$ for $n \in \{m+1, m+2, \ldots, 2m\}$.
Finally,
$$
\phi(2m+1) - \phi(2m+2)
= \sum_{\ell=0}^m { {m}\choose{\ell}} {{\ell}\choose{\ell+(2m+1-1-\ell)-m} }
= \sum_{\ell=0}^m { {m}\choose{\ell}} {{\ell}\choose{m} } = 1,
$$
so that
$\phi(n) > \phi(n+1)$ for $n \in \{m+1, m+2, \ldots, 2m+1\}$.
This finishes the proof.
\qed

\section{Depth}[sectdepth]%

The depth of quotients of powers of $\BaHH(m,r,s)$ depend on $s$,
so in this section we return to arbitrary $s$.

\lemma[lmfordepth]%
Set $B = \BaHH(m,r,s)$.
Let $w = a_1^{e_n} a_1^4 \cdots a_r^4 \prod_{i,j} x_{i,j}$
be defined as in the proof of \ref{thmnochigherpower}.
If $n\ge2$, then $w \not \in B^n$
and $w$ multiplies
$(a_j, x_{i,j} : i \in [m], j \in [r])$
(but not $c_1, \ldots, c_s$)
into $(B_0 + X)^n$.

Let $u_1, \ldots, u_s$ be linear forms with $u_j$
of the form $c_j$ minus a linear combination $d_j$
in the variables $x_{i,j'}$ as $i, j'$ vary in $[m]$ and $[r]$,
respectively.  Then
$w \not \in B^n + (u_1, \ldots, u_s)$
and
$w (a_j, x_{i,j} : i \in [m], j \in [r])$
$\in B^n + (u_1, \ldots, u_s)$.
\endb

\proof
The first paragraph is an immediate consequence of the proof of
\ref{thmnochigherpower}.

For the second paragraph,
it is still the case that $w$ multiplies
the $a_j$, $x_{i,j}$ (but not $c_1, \ldots, c_s$)
into the ideal $C = B^n + (u_1, \ldots, u_s)$.
It remains to prove that $w \not \in C$.
By \ref{lmpoweregen} we can rewrite $C$ as
$(c_1, \ldots, c_s) a_1^4 \cdots a_r^4 X^{n-1}
+ (B_0+X)^n + (u_1, \ldots, u_s)$
$= (d_1, \ldots, d_s) a_1^4 \cdots a_r^4 X^{n-1}
+ (B_0+X)^n + (u_1, \ldots, u_s)$.
Without restriction,
we can switch to the polynomial ring where $u_1, \ldots u_s$ are variables
and, for $1\le j\le s$, $c_j$ are the linear forms in $u_j$ and $d_j$.
Since $u_1, \ldots u_s$
do not appear in $w$ or in any minimal generating set of
$(d_1, \ldots, d_s) a_1^4 \cdots a_r^4 X^{n-1} + (B_0+X)^n$,
it follows that if $w$ is in $C$,
then
$w \in (d_1, \ldots, d_s) a_1^4 \cdots a_r^4 X^{n-1} + (B_0+X)^n$,
so that $w$ multiplies $c_1, \ldots, c_s$ into $B^n$,
which is a contradiction.
\qed

\thm[thmdepth2]%
For any positive integers $r, m, s, e$ with $r \ge 2$,
there exists an ideal $I$ in a polynomial ring $R$
such that for all positive integers~$n$,
$$
\depth\left({A \over I^n}\right)
= \cases{
e-1, & if $n = ru + 1$ with $u = 0, \ldots, m$; \cr
e, & if $n \le rm+1$ and $n \not \equiv 1 \mod r$; \cr
s+e-1, & otherwise, i.e., if $n > mr+1$. \cr
}
$$
In particular,
the depth function $n \mapsto \depth(A/I^n)$
has $m+1$ local minima,
it is periodic of period $r$ when restricted to the domain $[1, r(m+1)-1]$,
and it is constant afterwards.
\endb

\proof
Set $B = \BaHH(m,r,s)$ in the ambient polynomial ring~$R$.
Let $A$ be the polynomial ring obtained from $R$
by replacing the variable $c_1$ with variables $d_1, \ldots, d_e$.
Let $\varphi : R \to A$ be the algebra homomorphism
taking $c_1$ to the product $d_1 \cdots d_e$
and all other variables to themselves.
Let $I = \varphi(B) A$.
Since $\varphi$ is a free and hence a flat map by \cite[Theorem 1.2]{KS}
(such maps are called splittings there),
we have that $\varphi$ takes a free resolution of $R/B^n$
to a free resolution of $A/I^n$,
and it preserves minimality of the resolution.
Thus the projective dimensions of $R/B^n$ and $A/I^n$ are the same,
and by the Auslander-Buchsbaum formula,
$$
\eqalignno{
\depth(A/I^n)
&= \dim(A) - \projdim(A/I^n) \cr
&= \dim(R) + e - 1 - \projdim(R/B^n) \cr
&= \dim(R) + e - 1 - (\dim(R) - \depth(R/B^n)) \cr
&= \depth(R/B^n) + e - 1. \cr
}
$$
So it suffices to prove that
$$
\depth\!\left({R \over B^n}\right)
= \cases{
0, & if $n = ru + 1$ with $u = 0, \ldots, m$; \cr
1, & if $n \le rm+1$ and $n \not \equiv 1 \mod r$; \cr
s, & otherwise, i.e., if $n > mr+1$. \cr
}
$$
By \refs{thmspreadassoc} and \refn{thmcfullhalf},
the maximal ideal of $R$ is associated to $R/B^n$ exactly for the $n$
of the form $ru+1$ with $u = 0, \ldots, m$.
Thus the depth of $R/B^n$ equals $0$ exactly for all such $n$.

So we may assume that either $n \not \equiv 1 \mod r$ or that $n > mr + 1$.

No $c_1, \ldots, c_s$ appear in any generator of a minimal generating set of
$B^n : c_1 = (B + (a_1^4 \cdots a_r^4))^n$.
This means that $\depth(R/(B^n : c_1)) \ge s$.
By \ref{thmnochigherpower} (and \ref{thmspreadassoc})
we then have
$\depth(R/(B^n : c_1)) = s$.
Also,
$B^n + (c_1) = B(m,r,s-1)^n + (c_1)$
with $B(m,r,s-1)$ defined using variables $a_j, x_{i,j}$
and $c_2, \ldots, c_s$ only.
When $s = 1$,
by \ref{lmfordepth},
$\depth(R/(B^n + (c_1)))$ is $0$ for all~$n$,
and for $s \ge 2$,
by induction on $s$,
$\depth(R/(B^n + (c_1)))$ is $1$ or $s-1$, depending on $n$.
We will use the short exact sequence
$$
0 \longrightarrow {R \over B^n : c_1}
\longrightarrow {R \over B^n}
\longrightarrow {R \over B^n + (c_1)}
\longrightarrow 0.
\eqno{(*)}
$$
This short exact sequence induces a long exact sequence on $\Ext_R(R/M, \_)$,
where $M$ is the maximal homogeneous ideal of~$R$.
We use the fact that for any finitely generated $R$-module~$U$,
$\depth(U) = \min\{\ell: \Ext^\ell_R(R/M,U) \not = 0\}$.

Let $\ell = \depth(R/(B^n + (c_1))$.
By induction on $s$,
we have that $\ell = 0$ if $s = 1$,
and otherwise that $\ell = 1$ if $n \le mr+1$ (and $n \not \equiv 1 \mod r$),
and $\ell = s-1$ otherwise,
and so since the depth of $R/(B^n : c_1) = s$,
the relevant part of the long exact sequence equals:
$$
\cdots
\rightarrow
0
\rightarrow
\Ext_R^\ell \!\biggl({R\over M}, {R \over B^n} \biggr)
\rightarrow
\Ext_R^\ell \!\biggl({R\over M}, {R \over B^n + (c_1)} \biggr)
\rightarrow
\Ext_R^{\ell+1} \!\biggl({R\over M}, {R \over B^n: c_1} \biggr)
\rightarrow \cdots.
$$
We need to establish that
$\Ext_R^\ell\bigl({R\over M}, {R \over B^n} \bigr)$
is non-zero if $n \le mr+1$ (and $n \not \equiv 1 \mod r$)
and is zero if $n > mr+1$.

First let $n > mr + 1$.
We need to show that the depth of $R/B^n$ is~$s$.
By the long exact sequence we first prove that
$\Ext_R^{s-1} \!\left({R\over M}, {R \over B^n} \right)$ is zero,
i.e., that
$$
\Ext_R^{s-1} \!\left({R\over M}, {R \over B^n + (c_1)} \right)
\longrightarrow
\Ext_R^s \!\left({R\over M}, {R \over B^n: c_1} \right)
$$
is injective.
By faithful flatness we may assume that the base field is infinite.
By prime avoidance we can find linear forms $u_2, \ldots, u_s, u_1$
that form a regular sequence modulo $B^n : c_1$
and for which $u_2, \ldots, u_s$ is a regular sequence modulo $B^n + (c_1)$.
Since $a_1, \ldots, a_r$ are in the radical of $B$,
we may assume that the $u_i$
are forms in the variables $c_j$ and $x_{i,j}$ only.
We claim that for $\ell = 2, \ldots, s$, we may take $u_\ell = c_\ell - d_\ell$,
where $d_\ell = \sum_{i,j} \alpha_{\ell,i,j} x_{i,j}$
for some (generic) scalars $\alpha_{\ell,i,j}$.
Certainly any such $u_2, \ldots, u_s, u_1$ form a regular sequence
modulo $B^n : c_1$ since $c_1, \ldots, c_s$ do not appear
in any generators of this ideal.
Suppose that we have proved for some $\ell \in \{1, \ldots, s-1\}$
that $u_2, \ldots, u_\ell$ form a regular sequence modulo $B^n + (c_1)$.
Then
$$
B^n + (c_1, u_2, \ldots, u_\ell)
= \left(
(d_2, \ldots, d_\ell, c_{\ell+1}, \ldots, c_s) a_1^4 \cdots a_r^4
+ B_0 + X \right)^n + (c_1, u_2, \ldots, u_\ell),
$$
and since $(a_1^4 \cdots a_r^4)^2 \in B_0^2$,
by \ref{lmpoweregen},
modulo the variables $c_1, u_2, \ldots u_\ell$,
$$
B^n
=
(c_{\ell+1}, \ldots, c_s) a_1^4 \cdots a_r^4 X^{n-1} +
\left(
(d_2, \ldots, d_\ell) a_1^4 \cdots a_r^4 + B_0 + X \right)^n.
$$
Thus by \ref{lmoursuffices},
each associated prime of $B^n + (c_1, u_2, \ldots, u_\ell)$
either contains all $c_{\ell+1}, \ldots, c_s$ or it contains none of them.
Thus $c_{\ell+1} - d_{\ell+1}$
for sufficiently general $\alpha_{\ell+1,i,j}$ is a non-zerodivisor
modulo $B^n + (u_2, \ldots, u_\ell)$.
This proves the stated forms of $u_2, \ldots, u_s$
and we may also take $u_1 = c_1$.

By a theorem of Rees (see \cite[Lemma 2 (i)]{Matsu}),
due to natural isomorphisms,
it suffices to prove that the natural map
$$
\Hom_R \!\left({R\over M}, {R \over B^n + (c_1,u_2,\ldots, u_s)} \right)
\longrightarrow
\Hom_R \!\left({R\over M}, {R \over (B^n: c_1) + (c_1, u_2, \ldots, u_s)} \right)
$$
is injective.
In other words,
we need to show that the natural map
${L_1 : M \over L_1} \rightarrow
{L_2 : M \over L_2}$ 
is injective,
where
$L_1 = B^n + (c_1, u_2, \ldots, u_s)$
and
$L_2 = (B^n : c_1) + (c_1, u_2, \ldots, u_s)$.
Let $w \in (L_1 : M) \cap L_2$.
We have to prove that $w \in L_1$.

By subtracting elements of $L_1$ from~$w$,
by \ref{lmpoweregen},
$w \in a_1^4 \cdots a_r^4 X^{n-1}$,
and since $(a_1, \cdots, a_r, c_1, \ldots, c_s) a_1^4 \cdots a_r^4 X^{n-1}
 \subseteq B^n \subseteq L_1$,
we may assume that
$$
w = \sum_{\nu} e_\nu
a_1^4 \cdots a_r^4 \biggl(\prod_{i,j} x_{i,j}^{v_{\nu,i,j}}\biggr)
\biggl(\prod_{i,j} h_{i,j}^{u_{\nu,i,j}}\biggr),
$$
where for all $\nu$,
$e_\nu \in k$ and
$\sum_{i,j} u_{\nu,i,j} = n-1$.
Since $L_1, L_2$ are not monomial ideals, $w$ need not be a monomial.
Nevertheless,
we claim that each $x_{i,j}$
multiplies each summand in $w$ into $L_1$.
Proof of the claim:
Fix $(i,j)$.
We know that $x_{i,j} w \in L_1$.
This means that in at least one monomial summand $w_0$ of $w$,
$x_{i,j}$ must be incorporated into some new factor of $B$
while at the same time possibly breaking up some of the $n-1$ factors that are
generators of $B$.
Then by \ref{coraxn}~(1),
$x_{i,j} w_0 \in B^n \subseteq L_1$.
Hence $x_{i,j} (w - w_0)$
is also in $L_1$, and $w-w_0$ has fewer summands in it,
which proves the claim for $(i,j)$
by induction on the number of monomial summands in $w$.
This proves that every monomial appearing in~$w$
is multiplied by $(x_{i,j}: i \in [m], j \in [r])$ into~$B^n$
Thus by \ref{coraxn}~(2),
each monomial appearing in $w$ is in $L_1$,
so that $w \in L_1$.
This proves that $\Ext^{s-1}_R(R/M, R/B^n) = 0$,
which means that the depth of $R/B^n$ is at least~$s$.
By the same reasoning as before,
there exists a regular sequence $u_1, \ldots, u_s$ on $R/B^n$
of the form $u_j = c_j - d_j$
for some generic linear combinations $d_1, \ldots, d_s$ in the $x_{i,j}$.
Consider the element $w = a_1^{e_n} a_1^4 \cdots a_r^4 \prod_{i,j} x_{i,j}$,
with $e_n$ as defined in the proof of \ref{thmnochigherpower}.
By \ref{lmfordepth},
$Mw \in B^n + (u_1, \ldots, u_s)$
and $w \not \in B^n + (u_1, \ldots, u_s)$,
so that the depth of $R/B^n$ is exactly~$s$.

Finally, let $n \le mr+1$ and $n \not \equiv 1 \mod r$.
We need to prove that the connecting homomorphism in the displayed long
exact sequence is not injective.
If $s\ge3$,
then $\Ext_R^{\ell+1}\!\biggl({R\over M}, {R \over B^n: c_1} \biggr) = 0$
and we are done.
So, let $s\le2$.
As in in the proof for $n > mr+1$
there exists a linear form $u_2 = c_2 - d_2 \in M$ that
is a non-zerodivisor modulo $B^n : c_1$ and modulo $B^n + (c_1)$,
and by the same theorem of Rees, due to natural isomorphisms,
it suffices to prove that the natural homomorphism
${L_1 : M \over L_1} \rightarrow {L_2 : M \over L_2}$ 
is not injective,
where $L_1 = B^n + (c_1, u_2)$
and $L_2 = (B_n : c_1) + (c_1, u_2)$.
Write $n-2 = ur + v$ for some non-negative integers $u, v$
with $v < r$.
Since $n \le mr+1$ and $n \not \equiv 1 \mod r$,
it follows that $u < m$ and $v \not \equiv 0 \mod r$.
Let $w_h = h_{1,1} \cdots h_{1,v} \prod_{i > m-u, j} h_{i,j} \in B^{n-2}$,
$w_x = \prod_{i \le m-u, j} x_{i,j}$,
$w_a = a_1^4 \cdots a_{r-1}^4 a_r^8 \in B$,
and $w = x_{1,1} w_h w_a w_x$.
Thus clearly
$w \in B^{n-2} (a_r^4 x_{1,r} x_{1,1}^2) (a_1^4 \cdots a_r^4) \in L_2$.
However, $w \not \in L_1$.
We next prove that $M w \subseteq L_1$:
$$
\eqalignno{
x_{i,j} w &\in {w_h \over h_{i,j-2} h_{i,j}} h_{i,j-1} (a_{j-2}^6) (a_j^6)
(h_{1,r}) \subseteq (B_0 + X)^n,
\hbox{ if $i > m-u$}; \cr
x_{i,j} w &\in (w_h) (a_{j-1}^4 x_{i,j-1} x_{i,j}^2) (a_r^6) \subseteq (B_0 + X)^n,
\hbox{ if $j \not = 1$ and $i \le m-u$}; \cr
x_{i,1} w &\in (w_h) (a_r^4 x_{i,r-1} x_{i,1}^2) (a_r^4 x_{1,r-1} x_{1,1}^2)
\subseteq (B_0 + X)^n,
\hbox{ if $2 \le i \le m-u$}; \cr
x_{1,1} w &\in
(c_2 a_1^4 \cdots a_r^4) w_h (a_r^4 x_{1,r} x_{1,1}^2) + (u_2) + \sum_{(i,j) \not = (1,1)} x_{i,j} w \in L_1; \cr
a_r w &\in (w_h) (a_r^5 a_1) (a_r^4 x_{1,r} x_{1,1}^2) \in L_1; \cr
a_j w &\in (w_h) (a_j^5 a_{j+1}) (a_r^4 x_{1,r} x_{1,1}^2) \in L_1, \hbox{if
$j \not = 1$}. \cr
}
$$
This proves that $w \in (L_1 : M) \cap L_2$ and $w \not \in L_1$,
which proves that the map
${L_1 : M \over L_1} \rightarrow {L_2 : M \over L_2}$ 
is not injective.
Thus the depth of $R/B^n$ is~$1$
if $n \le rm+1$ and $n \not \equiv 1 \mod r$.
\qed

\bigskip\bigskip
\leftline{\bf References}
\bigskip

\baselineskip=9.9pt
\parindent=3.6em
\setbox1=\hbox{[999]} 
\newdimen\labelwidth \labelwidth=\wd1 \advance\labelwidth by 2.5em
\ifnumberbibl\advance\labelwidth by -2em\fi

\thmno=0

\bitem{AH}
T. Ananyan and M. Hochster,
Small subalgebras of polynomial rings and Stillman's conjecture,
{\it J. Amer. Math. Soc.} {\bf 33} (2020), 291--309.

\bitem{BHH}
S. Bandari, J. Herzog, and T. Hibi,
Monomial ideals whose depth function
has any given number of strict local maxima,
{\it Ark. Mat.} {\bf 52} (2014), 11--19.

\bitem{Magma}
W. Bosma, J. Cannon, C. Playoust,
The Magma algebra system. I. The user language
{\it J. Symbolic Comput.} {\bf 24} (1997), 235--265.

\bitem{Brodmann}
M. Brodmann,
Asymptotic stability of Ass($M/I^nM$),
{\it Proc. Amer. Math. Soc.} {\bf 74} (1979), 16--18.

\bitem{GS}
D. Grayson and M. Stillman,
Macaulay2, a software system for research in algebraic geometry,
available at {\tt http://www.math.uiuc.edu/Macaulay2}.

\bitem{HNTT}
H. T. H\`a, H. D. Nguyen, N. V. Trung, and  T. N. Trung,,
Depth functions of powers of homogeneous ideals,
{\it Proc. Amer. Math. Soc.} {\bf 149} (2021), 1837--1844.

\bitem{HH}
J. Herzog and T. Hibi,
The depth of powers of an ideal,
{\it J. Algebra} {\bf 291} (2005), 534--550.

\bitem{KS}
J. Kim and I. Swanson,
Many associated primes of powers of prime ideals,
{\it Journal of Pure and Applied Algebra} {\bf 223} (2019), 4888--4900.

\bitem{Matsu}
H. Matsumura,
{\bkt Commutative Ring Theory},
Cambridge University Press, 1986.

\bitem{Sei74}
A. Seidenberg,
Constructions in algebra,
{\it Trans. Amer. Math. Soc.} {\bf 197} (1974), 273--313.

\bitem{WS}
S. J. Weinstein and I. Swanson,
Predicted decay ideals,
{\it Comm. in Algebra}, {\bf 48} (2020), 1089--1098.
Published online on 26 October 2019.
{\tt arXiv:1808.09030}.

\end